\documentclass[11pt,reqno]{amsart}
\usepackage{amsfonts}
\usepackage{mathrsfs}
\usepackage{graphicx} 
\usepackage{epsfig}
\usepackage{amsmath, amsthm}
\usepackage{amssymb}

\theoremstyle{plain}
\newtheorem{thm}{Theorem}[section]
\newtheorem{prop}[thm]{Proposition}
\newtheorem{lem}[thm]{Lemma}
\newtheorem{cor}[thm]{Corollary}

\theoremstyle{definition}

\newtheorem{rem}[thm]{Remark}
\newtheorem{defn}[thm]{Definition}
\newtheorem{eg}[thm]{Example}
\newtheorem{subtitle}[thm]{}
\newtheorem{ex}{Exercise}[section]
\numberwithin{equation}{section}

\def\d{\delta}
\def\D{\triangle}
\def\e{\epsilon}
\def\g{\gamma}

\def\K{\nabla}
\def\l{\lambda}

\def\n{\vert\/}

\def\li{\langle}
\def\ri{\rangle}
\def\n{|\/ }
\def\tr{{\rm tr}}
\def\bs{\bigskip}

\def\ni{\noindent}
\def\ti{\tilde}
\def\p{\partial}

\def\I{{\rm I\/}}

\newcommand{\beq}{\begin{equation}}
\newcommand{\eeq}{\end{equation}}
\newcommand{\beg}{\begin{eg}}
\newcommand{\eeg}{\end{eg}}
\newcommand{\bthm}{\begin{thm}}
\newcommand{\ethm}{\end{thm}}
\newcommand{\bprop}{\begin{prop}}
\newcommand{\eprop}{\end{prop}}
\newcommand{\bcor}{\begin{cor}}
\newcommand{\ecor}{\end{cor}}
\newcommand{\blem}{\begin{lem}}
\newcommand{\elem}{\end{lem}}
\newcommand{\bca}{\begin{cases}}
\newcommand{\eca}{\end{cases}}
\newcommand{\brem}{\begin{rem}}
\newcommand{\erem}{\end{rem}}
\newcommand{\bpm}{\begin{pmatrix}}
\newcommand{\epm}{\end{pmatrix}}
\newcommand{\bbm}{\begin{bmatrix}}
\newcommand{\ebm}{\end{bmatrix}}
\newcommand{\bvm}{\begin{vmatrix}}
\newcommand{\evm}{\end{vmatrix}}
\newcommand{\bdefn}{\begin{defn}}
\newcommand{\edefn}{\end{defn}}
\newcommand{\bsub}{\begin{subtitle}}
\newcommand{\esub}{\end{subtitle}}
\newcommand{\bex}{\begin{ex}}
\newcommand{\eex}{\end{ex}}
\newcommand{\ben}{\begin{enumerate}}
\newcommand{\een}{\end{enumerate}}

\date{}

\def\calB{\mathcal{B}}

\def\calF{\mathcal{F}}
\def\calG{\mathcal{G}}

\def\calL{\mathcal{L}}
\def\calM{\mathcal{M}}
\def\calN{\mathcal{N}}
\def\calO{\mathcal{O}}
\def\R{\mathbb{R}}

\def\calT{\mathcal{T}}

\def\cf{{\mathcal{F}}}

\def\cl{{\mathcal{L}}}
\def\cm{{\mathcal{M}}}

\def\co{{\mathcal{O}}}

\def\R{\mathbb{R}}

\def\C{\mathbb{C}}

\def\fk{\mathfrak{k}}

\def\det{{\rm det \/ }}

\def\mod{{\rm mod\,}}
\def\n{\, | \,}

\def\rd{{\rm \/ d\/}}

\def\sech{{\rm sech\/}}

\def\an1{A^{(1)}_{n-1}}
\def\Ker{{\rm Ker\/}}

\def\r0{\R^n \backslash \{0\}}

\def\bh{\backslash}

\begin{document}

\title[Affine Curve Flow]
{N-dimension Central Affine Curve Flows} \today

\author{Chuu-Lian Terng$^\dag$}\thanks{$^\dag$Research supported
in  part by NSF Grant DMS-1109342}
\address{Department of Mathematics\\
University of California at Irvine, Irvine, CA 92697-3875.  Email: cterng@math.uci.edu}
\author{Zhiwei Wu$^*$}\thanks{$^*$Research supported in part by NSF of China under Grant No. 11401327\/}
\address{Department of Mathematics\\ Ningbo University\\ Ningbo, Zhejiang, 315211, China. Email: wuzhiwei@nbu.edu.cn}


\maketitle

\section{Introduction}

The group $SL(n, \R)$ acts on $\R^n \backslash \{0\}$
transitively by $g \cdot y = g y$ for $g \in SL(n, \R)$ and $y \in \R^n$. Given a curve $\g$ in $\r0$, if $\det(\g, \g_s, \ldots, \g_s^{(n-1)})$ is positive then there is an orientation preserving parameter $x$ unique up to translation, $\frac{\rd x}{\rd s}= \det(\g, \g_s, \ldots, \g_s^{(n-1)})^{\frac{2}{n(n-1)}}$,  such that 
\beq\label{do1}
\det(\g, \g_x, \ldots, \g_x^{(n-1)})=1,
\eeq
where $\g_x^{(i)}= \frac{\rd^i \g}{\rd x^i}$.  Take $x$ derivative of \eqref{do1} to get
 $$\det(\g, \g_x, \ldots, \g_x^{(n-2)}, \g_x^{(n)})=0.$$  Hence we have
$$
\g_x^{(n)}=u_1\gamma+u_2\gamma_x+\cdots+u_{n-1}\g_x^{(n-2)},
$$
where
$u_i=\det(\g, \g_x, \ldots, \g_x^{(i-2)},  \g_x^{(n)}, \g_x^{(i)}, \ldots, \g_x^{(n-1)})$.  This parameter $x$ is called the {\it central affine arc-length parameter\/}, 
$$g=(\g, \ldots, \g_x^{(n-1)})$$ the {\it central affine moving frame\/}, and $u_i$ the {\it $i$-th central affine curvature\/} of $\g$ for $1\leq i\leq n-1$ (cf. \cite{CIM13}, \cite{UP95}).  
Note that 
$$
g_x=g(b+u),
$$
where $b=\sum_{i=1}^{n-1}e_{i+1, i}$ and $u=\sum_{i=1}^{n-1}u_ie_{in}$.  We also call $u$ the {\it central affine curvature\/} along $\g$.

Let $\I= S^1$ or $\R$, and 
\begin{align*}
&\calM_n(\I)=\{\g:\I\to \r0 \n \det(\g,\g_x, \ldots, \g_x^{(n-1)})=1  \},\\
& V_n=\oplus_{i=1}^{n-1}\R e_{in} \subset sl(n,\R).
\end{align*}
Let $\Psi: \calM_n(\I)\to C^\infty(\I, V_n)$ be the {\it central affine curvature map\/} defined by 
$$
\Psi(\g)= u= \sum_{i=1}^{n-1} u_i e_{in},
$$
where $u_1, \ldots, u_{n-1}$ are the central affine curvatures along $\g$.  

It follows from the Existence and Uniqueness for Ordinary Differential Equations that $\Psi$ induces a bijection from the orbit space $\frac{\calM_n(\R)}{SL(n,\R)}$ to $C^\infty(\R, V_n)$ and $ u_1, \ldots, u_{n-1}$ form a complete set of local differential invariants for curves in $\r0$ under the group $SL(n,\R)$.  For example, curves with zero central affine curvatures in $\R^n\backslash \{0\}$  are of the form
$$\gamma(x)=c\left(1,x,\frac{x^2}{2}, \cdots, \frac{x^{n-1}}{(n-1)!}\right)^t$$
for some $c\in SL(n,\R)$. 

Let $y\in C^\infty(\R, \R^n)$. We say $\eta$ is a {\it differential polynomial of order $k$ in $y$} if $\eta$ is a polynomial in $y, y_x, \ldots, y_x^{(k)}$. 

A vector field $X:\cm_n(\R)\to C^\infty(\R, \R^n)$ is tangent to $\cm_n(\R)$ if and only if 
\beq\label{tf}
\sum_{i=0}^{n-1} \det(\g, \ldots, \g_x^{(i-1)}, (X(\g))_x^{(i)}, \ldots, \g_x^{(n-1)})=0.
\eeq
We show in section \ref{bs} that there exists a differential polynomial $\phi_n$ in $u$, $\xi_1$, \ldots, $\xi_{n-1}$ such that the vector field 
$\sum_{i=0}^{n-1}\xi_i \g_x^{(i)}$ is tangent to $\cm_n(\R)$ at $\g$ if and only if 
$$\xi_0= \phi_n(u, \xi_1, \ldots, \xi_{n-1}),$$ 
 where $u$ is the central affine curvature along $\g$. 
 So $T\calM_n(\R)_\g$ is identified as $C^\infty(\R, \R^{n-1})$.
 
A {\it central affine curve flow\/} is an evolution equation on $\calM_n(\R)$ of the form
\beq\label{age}
\g_t =X(\g)= \xi_0(u)\g + \xi_1(u)\g_x + \cdots + \xi_{n-1}(u)\g_x^{(n-1)},
\eeq
where $X(\g)$ is tangent to $\calM_n(\R)$ at $\g$ and $\xi_0, \ldots, \xi_{n-1}$ are differential polynomials in the central affine curvature $u(\cdot, t)$ of $\g(\cdot, t)$. So $\xi_0(u)= \phi_n(u, \xi_1(u), \ldots, \xi_{n-1}(u))$.

 It is easy to see that a central affine curve flow is invariant under the action of $SL(n,\R)$ on $\r0$ and translations in the $(x,t)$-plane. In other words, if $\g(x,t)$ is a solution of \eqref{age}, then so is $\ti\g(x,t)= c\g(x+r_1, t+r_2)$, where $c\in SL(n,\R)$ and $r_1, r_2\in \R$ are constants.   
 
 Let $n\geq 3$.  Equation \eqref{tf} implies that $X(\g)= y_1(\g) \g+ \g_{xx}$ is tangent to $\calM_n(\R)$ at $\g$ if and only if $y_1(\g)= -\frac{2}{n}u_{n-1}$. So
 \beq\label{ncf2}
\gamma_t=-\frac{2}{n}u_{n-1}\gamma+\gamma_{xx},
\eeq
is one of the simplest central affine curve flows on $\calM_n(\R)$ for $n \geq 3$,
where $u_{n-1}(\cdot, t)$  is the $(n-1)$-th central affine curvature for $\gamma(\cdot, t)$.
This curve flow turns out to be integrable. In fact, if $\g$ is a solution of the central affine curve flow \eqref{ncf2} then its central affine curvature $u(\cdot, t)$ is a solution of the second flow of the Drinfeld-Sokolov $\an1$-KdV hierarchy constructed in \cite{DS84}. The $\an1$-KdV hierarchy is the same as the Gelfand-Dickey (GD$_n$) hierarchy on the space of n-th order differential operators on the line. So the central affine curve flow \eqref{ncf2} can be viewed as a simple and natural geometric interpretation of the second GD$_n$ flow.

Next we state some of our main results:

\ben
\item[(a)] We construct a sequence of commuting higher order central affine curve flows on $\calM_n(\R)$ such that the second flow is \eqref{ncf2}. For example, the third central affine curve flow on $\calM_n(\R)$ ($n \neq 3$) is  
\beq\label{ki}
\gamma_t=\left(-\frac{3}{n}u_{n-2}+\frac{3(n-3)}{2n}(u_{n-1})_x\right)\gamma
- \frac{3}{n}u_{n-1}\gamma_x+\gamma_{xxx}.
\eeq
Note that when $n=2$, we have $\g_{xx}= u_1\g$. Then \eqref{ki} becomes
\beq\label{ki2}
\g_t= \frac{1}{4} (u_1)_x \g -\frac{1}{2} u_1\g_x.
\eeq
It was proved in \cite{UP95} that if $\g$ is a solution of \eqref{ki2} on $\calM_2(\R)$ then $u_1$ is a solution of the KdV equation.

\item[(b)] We prove that the central affine curvature map $\Psi$ gives a one to one correspondence between solutions of the $j$-th central affine curve flow modulo the action of $SL(n,\R)$ and solutions of the $j$-th $\an1$-KdV flow.  

\item[(c)] We use the solution of the Cauchy problem of the $\an1$-KdV hierarchy to solve the Cauchy problem for the central affine curve flow hierarchy with periodic initial data and also with initial data having rapidly decaying central affine curvatures.  

\item[(d)] We construct B\"acklund transformations, a permutability formula, and infinitely many families of explicit solutions of the central affine curve flows. 
\een

 We pull back the known sequence of Poisson structures for the $\an1$-KdV hierarchy via the central affine curvature map $\Psi$ to $\calM_n(S^1)$ for the central affine curve flow hierarchy. We also prove the following results in this paper:
 \ben
 \item The central affine curve flows \eqref{ncf2} and \eqref{ki} are the Hamiltonian equations for the functionals 
\begin{align*}
F_2(\g)&= \oint u_{n-2}(x) \rd x,\\
F_3(\g)&= \oint u_{n-3} + \frac{n-3}{2n} u_{n-1}^2 \rd x,
\end{align*}
with respect to the second Poisson structure respectively, where $u_i$ is the $i$-th central affine curvature along $\g$.  
 \item The second and third Poisson structures arise naturally from co-adjoint orbits. 
 \item We identify the kernels of these two Poisson operators.  
\item
 Since the $\R^{n-1}$-action on $\calM_n(S^1)$ generated by the first $(n-1)$ central affine curve flows commutes with the $SL(n,\R)$-action, the direct product $SL(n,\R)\times \R^{n-1}$ acts on $\calM_n(S^1)$. We show that  the second and third Poisson structures on $\calM_n(S^1)$  induces weak symplectic forms on the orbits space $\frac{\calM_n(S^1)}{SL(n,\R)}$  and $\frac{\calM_n(S^1)}{SL(n,\R)\times \R^{n-1}}$ respectively. 
 \een

When $n=3$, result (a) and one Poisson structure were obtained in \cite{CIM13} for \eqref{ncf2}.  When $n=2$, result (4) was proved in \cite{FK13}. 

This paper is organized as follows: In section \ref{bs}, we give a brief review of the $\an1$-KdV hierarchy and prove some properties of its Lax pair. We prove results (a)-(c) in section \ref{gg} and (d) in section \ref{op}.  We review the bi-Hamiltonian structure of the $\an1$-KdV hierarchy and compute the kernels of the bi-Hamiltonian structures in section \ref{go}. We write down the formula for the bi-Hamiltonian structures of the curve flow \eqref{ncf2} and prove results concerning Hamiltonian structures in the last section.

\bs

\section{The $A^{(1)}_{n-1}$-KdV hierarchy}\label{bs}

Drinfeld and Sokolov constructed a hierarchy of KdV type for each affine Kac-Moody algebra.  In this section, we 
\ben
\item[(i)]  give a brief review of the construction of the $\an1$-KdV hierarchy (cf. \cite{DS84}),
\item[(ii)] develop some properties of its Lax pair that are needed for the study of central affine curve flows and the bi-Hamiltonian structure,
\item[(iii)] identify $T\cm(\R)_\g$ as $C^\infty(\R, \R^{n-1})$.
\een

\medskip

We first set up some notations. 
 Let 
 \begin{align*}
 &\calB_n^+=\{y=(y_{ij})\in sl(n,\R)\n y_{ij}=0, i>j\},\\
 &\calN_n^+=\{y=(y_{ij})\in sl(n,\R)\n y_{ij}=0, i\geq j \},\\
 &\calT_n =\{y\in gl(n,\R)\n y_{ij}=0, i \neq j \},
 \end{align*}
denote the
subalgebras of upper triangular, strictly upper triangular matrices in $sl(n,\R)$ and diagonal matrices in $gl(n,\R)$ respectively, and $N_n^+$  the corresponding Lie  subgroup of $\calN_n^+$.

Let
$$ \calL(sl(n,\R))=\left\{\xi(\l)=\sum_{j\leq n_0} \xi_j \l^j \bigg| \xi_j\in sl(n,\R), n_0\,\, {\rm integer}\right\}.$$ 
For $\xi\in \calL(sl(n,\R))$, we use the following notation:
$$\xi_+=\sum_{j\geq 0} \xi_j \l^j, \quad \xi_-=\sum_{j < 0} \xi_j \l^j.$$ 
Let
$$
J=e_{1, n}\l + b, \quad b= \sum_{i=1}^{n-1}e_{i+1, i}.
$$

Given $u \in C^{\infty}(\R, \calB_n^+)$, a direct computation (cf. \cite{TU14b}) implies that there exists a unique $Y(u,\l) \in \cl(sl(n, \C))$ satisfying
\beq\label{ff}
\bca
[\p_x+J+u, Y(u,\l)]=0, \\
Y(u,\l)^n=\l \I_n.
\eca
\eeq
Moreover, the coefficients of the power series expansion of $Y(u,\l)$ in $\l$ are differential polynomials of $u$.

Given $j\not\equiv 0$ ($\mod n$), let $Y(u,\l)$ denote the solution of \eqref{ff}, and write
\beq\label{fm}
(Y(u,\l))^j=\sum_{-\infty}^{[\frac{j}{n}]+1}Y_{j,i}(u)\l^i.
\eeq
 It was known (cf. \cite{DS84}, \cite{TU14b}) that if $u\in C^\infty(\R, V_n)$ then  there is a unique $\zeta_j(u)\in \calN_n^+$ such that
\beq\label{fmb}
[\p_x+b+u, Y_{j,0}(u)- \zeta_j(u)]\in V_n=\oplus_{i=1}^{n-1} \R e_{in}
\eeq
 and entries of $\zeta_j(u)$ are differential polynomials in $u$.  
 Set 
\beq\label{by} 
Z_j(u,\l)=(Y(u,\l)^j)_+- \zeta_j(u)= \sum_{0 \leq i\leq [\frac{j}{n}]+1} Z_{j,i}(u) \l^i. 
\eeq
 The  {\it $j$-th $\an1$-KdV hierarchy\/} ($j\geq 0$ and $j\not\equiv 0$ $\mod(n)$) constructed by Drinfeld-Sokolov in \cite{DS84} is the following flow on $C^\infty(\R,V_n)$:
\beq\label{mgd0}
u_t=[\p_x+b+u, \, Z_{j,0}(u)].
\eeq

\bprop\label{fn} (\cite{DS84}) Given $u\in C^\infty(\R^2, V_n)$, then 
 the following statements are equivalent 
 \ben
\item[(i)] $u$ is a solution of the $j$-th $\an1$-KdV flow, \eqref{mgd0},
\item[(ii)] for all parameter $\l\in \C$, $u$ satifies
\begin{equation}\label{laxp}
[\p_x+ J+u, \p_{t_j}+ Z_j(u,\l)]=0,
\end{equation}
 (i.e. \eqref{laxp} is the Lax pair of \eqref{mgd0}),
\item[(iii)]
$[\p_x+ b+u, \, \p_t+ Z_{j,0}(u)]=0$,
which is the Lax pair \eqref{laxp} with parameter $\l=0$,
\item[(iv)] the following system is solvable for $E(x,t,\l)\in GL(n,\C)$:
\beq\label{id}
\bca E^{-1}E_x= J+u, \\
 E^{-1}E_t= Z_j(u,\l),\\
 \overline{E(x,t,\bar\l)} = E(x,t,\l).
 \eca \eeq
\een
\eprop

 We call a solution $E$ of \eqref{id} a {\it  frame\/} of the solution $u$ of the $j$-th $\an1$-KdV flow \eqref{mgd0}. Note that the third condition of \eqref{id} implies that $E(x,t,\l)\in GL(n,\R)$ for $\l\in \R$.

It follows from $\rd (\ln\det(E))= \tr(E^{-1}\rd E)$ and $\tr(J+u)=\tr(Z_j(u,\l))= 0$ that we have

\bcor\label{jv} If $E(x,t,\l)$ is a frame of a solution $u$ of the $j$-th $\an1$-KdV flow, then $\det(E(x,t,\l))$ is independent of $x, t$. 
\ecor

Next we derive some properties of $Z_{j,0}(u)$.  Note that $C=Z_{j,0}(u)$ satisfies 
\beq\label{jwb}
[\p_x+ b+ u, C]\in V_n.
\eeq
The following Theorem shows that if $C$ satisfies \eqref{jwb} then $C$ is determined by $\{C_{i1}\n 2\leq i\leq n\}$ or by $\{C_{ni}\n 1\leq i\leq n-1\}$. 

\bthm\label{jw} Let $u\in C^\infty(\R, V_n)$, and $C_j$ the $j$-th column of   $C=(C_{ij})\in C^\infty(\R, sl(n,\R))$.  If $C$ satisfies \eqref{jwb}, then 
\ben
\item[(i)]  $C_{j+1}= (C_j)_x+(b+u)C_j$ for $1\leq j\leq n-1$,
\item[(ii)] there exists a differential polynomial $\phi_n(u, C_{21}, \ldots, C_{n1})$ such that 
$$C_{11}= \phi_n(u, C_{21}, \ldots, C_{n1}),$$
\item[(iii)]  entries of $C$ are differential polynomials in $u, C_{21}, \ldots, C_{n1}$,
\item[(iv)] for $1\leq i\leq n-1$, we have
\beq\label{he} 
C_{ni} =C_{n-i+1, 1} + (i-1)(C_{n-i+2,1})_x+ \phi_i(C_{n1}, \ldots, C_{n-i+3,1}), 
\eeq
for some linear differential operators $\phi_i$ with differential polynomials in $u$ as coefficients,
\item[(v)] $C_{ij}$'s are differential polynomials in $u, C_{n1}, \ldots, C_{n, n-1}$. 
\een
\ethm

\begin{proof} (i) follows from $C_x+[b+u, C] \in V_n$ and a direct computation.

It follows from (i) and induction on $i$ that there exist differential polynomials $\psi_i(u, C_{21}, \ldots, C_{n1})$ such that 
$$C_{i+1, i+1}= C_{ii} + \psi_i(u, C_{21}, \ldots, C_{n1})$$
for $1\leq i\leq n-1$.  Since $\tr(C)=0$, we obtain (ii). 
Statements (iii)-(v) can be proved using (i) and induction.  
\end{proof}

\bcor\label{jwc}
Let $\g\in \calM_n(\R)$, $u=\sum_{i=1}^{n-1} u_i e_{in}$ the central affine curvature along $\g$, and $C=(C_{ij}): \R\to sl(n,\R)$ satisfying
\eqref{jwb}.  Then $\xi(\g)= \sum_{i=1}^{n} C_{i1} \g_x^{(i-1)}$ is tangent to $\calM_n(\R)$ at $\g$. 
\ecor
 
\begin{proof}
Let $g=(\g, \ldots, \g_x^{(n-1)})$ be the central affine moving frame along $\g$, $C_i$ the $i$-th column of $C$, and $\eta_i= gC_i$. Then $g^{-1}g_x= b+u$ and $\xi(\g)= gC_1=\eta_1$. Let $\rho=[\p_x+b+u, C]=C_x+[b+u, C]$. Then $\rho \in V_n$, and 
\beq\label{jx}
(gC)_x=g_xC+gC_x=gC(b+u)+g\rho.
\eeq
Since the first $(n-1)$ columns of $\rho$ are zero, \eqref{jx} implies that $\eta_{i+1}=(\eta_i)_x$ for $1 \leq i \leq n-1$. So we have 
$$(\xi(\g))_x^{(i-1)}= (\eta_1)_x^{(i-1)}= \eta_i = gC_i.$$
Hence $\det(\g, \ldots, \g_x^{(i-2)}, (\xi(\g))_x^{(i-1)}, \g_x^{(i)}, \ldots, \g_x^{(n-1)})= C_{ii}$.  Since $C$ is in $ sl(n,\R)$, we have $\sum_{i=1}^n C_{ii}=0$. So $\xi(\g)$ satisfies \eqref{tf}, i.e., 
$$\sum_{i=1}^n \det(\g, \ldots, \g_x^{(i-2)}, (\xi(\g))_x^{(i-1)}, \g_x^{(i)}, \ldots, \g_x^{(n-1)})=0.$$ Hence $\xi(\g)$ is tangent to $\calM(\R)$ at $\g$. 
\end{proof}

Henceforth we set
$$e_1=(1,0,\ldots, 0)^t\in \R^n.$$

Since $Z_{j,0}(u)$ satisfies \eqref{jwb}, we have the following.

\bcor\label{jwd} Given $\g\in \calM_n(\R)$, let $g$ and $u$ be the central affine moving frame and curvature of $\g$ respectively.  Then $\xi(\g)= gZ_{j,0}(u) e_1$ is a tangent vector field of $\calM_n(\R)$ and 
\beq\label{ncfj}
\g_t= gZ_{j,0}(u)e_1
\eeq
is a flow equation on $\calM_n(\R)$. 
\ecor

We call \eqref{ncfj}  the {\it $j$-th central affine curve flow\/}.

 A direct computation shows that the first column of $Z_{2,0}(u)$ is 
$$(-\frac{2}{n} u_{n-1}, 0, 1, 0, \ldots, 0)^t.$$  So the central affine curve flow \eqref{ncf2} can be written as 
$\g_t= gZ_{2,0}(u)e_1$.

\beg {\bf [Higher order central affine curve flows]}\par

We use \eqref{ff} and \eqref{fmb} to compute $Z_{j,0}(u)$ and see that  the fourth and the fifth central affine curve flows on $\calM_3(\R)$ are
\begin{align*}
&\gamma_t=-\frac{1}{9}(2u_2''-3u_1'-2u_2^2)\gamma+\frac{1}{3}(u_2'-u_1)\gamma_x
-\frac{u_2}{3}\gamma_{xx}, \\
&\gamma_t=\frac{1}{9}(-u_1''+u_1u_2)\gamma-\frac{1}{9}(u_2''-3u_1'+u_2^2)\gamma_x
+\frac{1}{3}(u_2'-2u_1)\gamma_{xx}.
\end{align*}
For $n\not=3$, the third central affine curve flow on $\calM_n(\R)$ is  the flow \eqref{ki}.
\eeg

Next we use Theorem \ref{jw} to define an operator $P_u$ that will be used in the construction of bi-Hamiltonian structure later. 

\bdefn\label{gp} Fix $u\in C^\infty(\R, V_n)$, let $V_n^t= \oplus_{i=1}^{n-1} \R e_{ni}$, and let $$P_u:C^\infty(\R, V_n^t)\to C^\infty(\R, sl(n,\R))$$ denote the map defined by $P_u(v)=C$, where $C$ is the unique $sl(n,\R)$-valued map satisfies $\pi_0(C)=v$ and $[\p_x+ b+u, C]\in C^\infty(\R, V_n)$, where $\pi_0$ is the {\it canonical projection\/} defined by 
$$\pi_0:sl(n,\C)\to V_n^t, \quad \pi_0(y)= \sum_{i=1}^{n-1} y_{ni} e_{ni}, \quad{\rm for} \quad y=(y_{ij}).$$ 
\edefn

It follows from Theorem \ref{jw} that the entries of $P_u(v)$ are differential polynomials in $u$ and $v$. 

The following is a consequence of Corollary \ref{jwc}. 

\bcor\label{sy}  Suppose $u, g$ are the central affine curvature and moving frame along $\g\in \calM_n(\R)$ respectively. Let $v\in C^\infty(\R, V_n^t)$. Then there exists $\d\g$ tangent to $\calM_n(\R)$ at $\g$ such that $P_u(v)= g^{-1}\d g$, where $\d g= (\d \g, (\d\g)_x, \ldots, (\d \g)_x^{(n-1)})$. 
\ecor

Since $[\p_x+b+u, Z_{j,0}(u)]\in V_n$, we have $Z_{j,0}(u)= P_u(\pi_0(Z_{j,0}(u))$.  It follows from $\zeta_j(u)\in \calN_n^+$ and $Z_{j,0}(u)= Y_{j,0}(u)-\zeta_j(u)$ that we have the following. 

\bcor\label{gpc} Let $\pi_0:sl(n,\R)\to V_n^t$ be the canonical projection,  $Y(u,\l)$ the solution of \eqref{ff}, and $Y_{j,0}(u), Z_{j,0}(u,\l)= Y_{j,0}(u,\l) -\zeta_j(u)$ as in \eqref{fm} and \eqref{by}.  Then
 $Z_{j,0}(u)=P_u(\pi_0(Z_{j,0}(u)))= P_u(\pi_0(Y_{j,0}(u)))$, and the $j$-th $\an1$-KdV flow \eqref{mgd0} on $C^\infty(\R, V_n)$ can be written as
 \beq\label{mgd}
u_{t_j}= [\p_x + b+u, P_u(\pi_0( Y_{j,0}(u)))].
\eeq
\ecor

Next we identify $T\cm(\R)_\g$ as $C^\infty(\R, \R^{n-1})$. 

\bcor\label{jwe} The vector field   $\sum_{i=0}^{n-1} \xi_i \g_x^{(i)}$ is tangent to $\cm_n(\R)$ at $\g\in \cm_n(\R)$ if and only if $\xi_0= \phi_n(u, \xi_1, \ldots, \xi_{n-1})$, where $u$ is the central affine curvature along $\g$ and $\phi_n$ is the differential polynomial in Theorem \ref{jw} (ii). 
\ecor

\begin{proof} 
 Given $\g\in \cm_n(\R)$, let $y:(-\e, \e)\to \cm_n(\R)$ with $y(0)=\g$.  Let $g(\cdot, s)$ and $u(\cdot, s)$  denote the central affine moving frame and curvature along $y(s)$ for each $s\in (-\e, \e)$. Let $\d\g= \frac{\rd y}{\rd s}|_{s=0}$, and $\d g= (\d \g, (\d\g)_x, \ldots, (\d\g)_x^{(n-1)})$.  Take $s$ derivative of $g^{-1}g_x=b+ u$ to see that $[\p_x+b+u, g^{-1}\d g]\in V_n$, i.e., $C=(C_{ij}):=g^{-1}\d g$ satisfies \eqref{jwb}. By Theorem \ref{jw}, we have $\d\g=gCe_1= \sum_{i=1}^n C_{i1} \g_x^{(i-1)}$ and $C_{11}= \phi_n(u, C_{21}, \ldots, C_{n1})$.  
 \end{proof}
 
 \beg
The proof of Theorem \ref{jw} gives an algorithm to compute the formula for $\phi_n$. 
 For example, we get $\sum_{i=0}^{n-1}\xi_i \g_x^{(i)}$ is tangent to $T\cm_n(\R)$ at $\g$ if 
\begin{align*}
& \xi_0= -\frac{1}{2} \xi_1', \quad {\rm for\, \,} n=2,\\
&\xi_0= -\frac{1}{3} (\xi_2''+ 3\xi_1 + 2u_2\xi_2), \quad  {\rm for\, \,} n=3,\\
& \xi_0= -\frac{1}{4} (\xi_3'''+ 4\xi_2'' + 6 \xi_1' + 3u_3' \xi_3+ 5 u_3\xi_3' + 2u_2\xi_3 + 2 u_3\xi_2), \quad  {\rm for\, \,} n=4.
\end{align*}
For $n=5$, we have
\begin{align*}
&\xi_0=-\frac{1}{5}(\xi_4^{(4)}+5\xi_3^{(3)}+10\xi_2''+10\xi_1'+6u_3'\xi_4+9u_3\xi_4'+4(u_4\xi_4)''+3(u_4\xi_4')' \\
&\qquad+2u_4\xi_4''+3(u_4\xi_3)'+4u_4\xi_3'+4u_2\xi_4+3u_3\xi_3+2u_4\xi_2+2u_4^2\xi_4) .
\end{align*}
Here we use $y'= y_x$, $y^{''}= y_x^{(2)}$, $y^{(3)}= y_x^{(3)}$ etc.
 For general $n$, we have 
$$\xi_0= -\frac{1}{n} ((\xi_{n-1})_x^{(n-1)} + \cdots).$$
\eeg

\bigskip

\section{Central affine curve flows and the $\an 1$-KdV hierarchy}\label{gg}

In this section, we prove the following results:
\ben
 \item  The affine curvature map $\Psi(\g)= u$ gives a one to one correspondence between the space of solutions of the $j$-th central affine curve flow \eqref{ncfj} on $\r0$ modulo $SL(n,\R)$ and the space of solutions of the $j$-th $\an1$-KdV flow,
\beq\label{dsj}
u_t= [\p_x+ b+ u, Z_{j,0}(u)].
\eeq
When $n=3$, this results was obtained in \cite{CIM13}. 
\item  We use solutions of the Cauchy problem for the j-th $\an1$-KdV flow to solve the Cauchy problem for \eqref{ncfj} with periodic initial data or initial data having rapidly decaying central affine curvatures.
\een

\bthm\label{cfagdnc}
Let $u=\sum_{i=1}^{n-1}u_ie_{in}$ be a solution of the $j$-th $\an1$-KdV flow \eqref{dsj}, and $c_0\in SL(n,\R)$ a constant. Let $g:\R^2\to SL(n,\R)$ denote the solution of  
\beq\label{idb}
 g^{-1}g_x= b+u,\quad g^{-1}g_t= Z_{j,0}(u),
\eeq
with $ g(0,0)= c_0$. Then 
$\gamma(x, t):= g(x,t)e_1$  is a
solution of the $j$-th central affine curve flow \eqref{ncfj} with central affine curvature $u(x, t)$.
\ethm

\begin{proof} Note that $g^{-1}g_x=b+u$ implies $g=(\g, \g_x, \ldots, \g_x^{(n-1)})$ and $\g_x^{(n)}= u_1 \g+ \ldots + u_{n-1} \g_x^{(n-2)}$.  So $\g(\cdot, t)\in \calM_n(\R)$ and $u_1, \ldots, u_{n-1}$ are the central affine curvatures of $\g(\cdot, t)$.   
Since $g_t= gZ_{j,0}(u)$,  we get
$\gamma_t=g_te_1=gZ_{j,
0}(u)e_1$, which is the $j$-th central affine curve flow. 
\end{proof}

The converse is also true.

\bthm\label{cfagdn} Let $\gamma$ be a solution of \eqref{ncfj} on $\calM_n(\R)$, and  $u(\cdot, t)= \sum_{i=1}^{n-1}u_i(\cdot, t) e_{in}$ the central affine curvature along $\g(\cdot, t)$. Then $u$ is a solution of the $j$-th $\an1$-KdV flow \eqref{dsj}.
\ethm

\begin{proof}
Let $g(\cdot, t)=(\gamma, \gamma_x, \cdots, \gamma_x^{(n-1)})(\cdot, t)$ be the central affine moving
frame for $\gamma(\cdot, t) \in \calM_n(\R)$.  Then $g^{-1}g_x= b+u$.
Next we compute $g_t$. Let $\eta_i$ denote the $i$-th column of $gZ_{j,0}(u)$. Then \eqref{ncfj} is $\g_t=\eta_1$.  Since $[\p_x+b+u, Z_{j,0}(u)]\in V_n$, Theorem \ref{jw} implies that $\eta_i= (\eta_1)_x^{(i-1)}$. A direct computation implies that  $(\g_x^{(i-1)})_t= (\g_t)_x^{(i-1)}= (\eta_1)_x^{(i-1)}$ for $2\leq i\leq n$.  This proves that $g_t= gZ_{j,0}(u)$.  It follows from Proposition \ref{fn} that $u$ is a solution of \eqref{dsj}.
\end{proof}

\bcor Let $\Psi$ denote the central affine curvature map, and $\g_1, \g_2$ solutions of \eqref{ncfj} on $\calM_n(\R)$. Then
\ben
\item $\Psi(\g_1(\cdot, t))= \Psi(\g_2(\cdot, t))$ if and only if there is a constant $c_0$ in $SL(n,\R)$ such that $\g_2=c_0\g_1$,
\item $\Psi$ induces a bijection between the space of solutions of \eqref{ncfj} modulo $SL(n,\R)$ and the space of solutions of \eqref{dsj}.
\een 
\ecor

Next we discuss the Cauchy problem for the $j$-th central affine curve flow \eqref{ncfj}. 
The Cauchy problem for the $\an 1$-KdV hierarchy \eqref{mgd} is
solved for an open dense subset of rapidly decaying smooth initial data using the method of inverse scattering (cf. \cite{BDT}). As a consequence,  we get the solution for the Cauchy problem for the curve flow \eqref{ncfj}:

\bthm {\rm [Cauchy problem with rapidly decaying affine curvatures]} \ 

\ni Let $j\geq 1$ and $j\not\equiv 0 \, (\mod n)$. Given $\gamma_0\in \calM_n(\R)$ with rapidly decaying central affine curvatures $u_1^0, \cdots, u_{n-1}^0$, let $g_0$ be the central affine moving frame along $\g_0$. 
Suppose $u=\sum_{i=1}^{n-1} u_i e_{in}$ is the solution of the $j$-th flow \eqref{dsj} in the $\an 1$-KdV hierarchy  with $u(x,0)=\sum_{i=1}^{n-1} u_i^0(x)e_{in}$. Let $g(x, t):\R^2\to SL(n,\R)$ be the solution of \eqref{idb} with initial data $g(0, 0)=g_0(0)$.  Then $\gamma=ge_1$ is the solution of the $j$-th central affine curve flow \eqref{ncfj} with $\gamma(x, 0)=\gamma_0(x)$. Moreover, the central affine curvatures of $\g(\cdot, t)$ are also rapidly decaying.
\ethm

Finally, we use the solution of Cauchy problem of the second $\an1$-KdV flow with periodic initial data to solve the Cauchy problem for the curve flow \eqref{ncf2} with periodic initial data.  By Theorem \ref{cfagdnc}, we only need to solve the period problem of \eqref{idb}. In fact, we have the following.

\bthm \label{c} {\rm [Cauchy Problem with periodic initial data]\/} \ 

\ni Let $\gamma_0\in \calM_n(S^1)$, and $u_1^0, \ldots, u_{n-1}^0$ the central affine curvatures of $\g_0$. Suppose
$u=\sum_{i=1}^{n-1} u_i e_{in}$ is the solution of the periodic Cauchy problem of \eqref{dsj} with  initial data $u(x,0)=\sum_{i=1}^{n-1} u_i^0e_{in}$.
Let $g:\R^2\to SL(n,\R)$ be the solution of \eqref{idb} with initial data $c_0=g_0(0)$, where $g_0$ is the central affine frame along $\g_0$. Then $\g= ge_1$ is a solution of \eqref{ncfj} with initial data $\g(x,0)= \g_0(x)$. Moreover, $\g(x,t)$ is periodic in $x$ and $\{u_i(\cdot, t), 1\leq i \leq n-1\}$ are the central affine curvatures for $\g(\cdot, t)$. 
\ethm 

\begin{proof} Note that both $g_0$ and $g(\cdot, 0)$ satisfy the same ordinary differential equation, $g^{-1}g_x= b+u(x, 0)$, and have the same initial data. So the uniqueness of ordinary differential equations implies that $g(x,0)=g_0(x)$. 
It follows from Theorem \ref{cfagdnc} that $\gamma(x,t)=g(x, t)e_1$ is a
solution of the curve flow \eqref{ncfj}. Moreover, 
$\gamma(x,0)= g(x,0)e_1= \gamma_0(x)$. It remains to prove that $\gamma$ is
periodic in $x$.

Since $\gamma_0$ is periodic with period $2\pi$, $g_0$ and $u_0$ are periodic in $x$ with period $2\pi$. Since $u(x, t)$ is periodic in $x$, so is $Z_{j, 0}(u)$. It suffices to prove 
$$y(t)= g(2\pi, t)- g(0, t)$$
is identically zero. To do this, we calculate
\begin{align*}
y_t&=g_t(2\pi, t)-g_t(0, t) \\
&=(gZ_{j,0}(u))(2\pi, t)-(gZ_{j,0}(u))(0, t)=(g(2\pi, t)-g(0, t))Z_{j,0}(u(0, t))\\
&=y(t)Z_{j,0}(u(0, t)).
\end{align*}
Since $g_0$ is periodic in $x$ with period $2\pi$, $y(0)=g(2\pi, 0)-g(0, 0)=0$. Note that $Z_{j,0}(u(0, t))$ is given and $0$ is the solution of $y_t= yZ_{j,0}(u(0, t))$ with the same initial condition $y(0)=0$. So it follows from the uniqueness of ordinary differential equations that $y$ is identically zero. 
\end{proof}

\bs

\section{B\"acklund transformations}\label{op}

In this section, we use B\"acklund transformations (BTs) for the $\an1$-KdV hierarchy given in \cite{TWc} to  construct BTs for the $j$-th central affine curve flow \eqref{ncfj}, a permutability formula, and infinitely many explicit rational and soliton solutions of the curve flow \eqref{ncfj}. 

We first summarize results concerning BTs of the $\an1$-KdV hierarchy obtained in \cite{TWc}.
Let $E$ and $\ti E$ be frames of solutions $u=\sum_{i=1}^{n-1} u_i e_{in}$ and $\ti u=\sum_{i=1}^n \ti u_i e_{in}$ of the $j$-th $\an1$-KdV flow \eqref{dsj} respectively. Suppose $\ti E= Ef$ and the first column of $f(x,t,\l)$ is $(h,1,0,\ldots,)^t$ for some $h\in C^\infty(\R^2, \R)$.  Then there exist differential polynomials $s_i(u,h)$, $\eta_{n,j}(u,h)$, and $\xi_n(u,h)$ such that
\ben
\item[(i)]
\beq\label{nm1}
\ti u_i= u_i + s_i(u,h), \quad 1\leq i\leq n-1,
\eeq
\item[(ii)]  $f$ is determined by $h$, in fact, $f=J+ h\I_n + N(u,h)$, where $N$ is strictly upper triangular and the $ij$-th entry of $N(u,h)$ is
$$
\bca N_{ij}(u,h)= C_{j-1, i-1} h_x^{(j-i)}, &1 \leq i < j < n,\\
N_{in}=u_i + s_i(u,h) + C_{n-1, i-1} h_x^{(n-i)}, & 1\leq i\leq n-1,\eca
$$
and $C_{j,i}=\frac{j!}{i!(j-i)!}$. Henceforth we use $f_{u,h}$ to denote such $f$.
\item[(iii)]  $\det(f_{u,h})= (-1)^{n-1}(\l + h_x^{(n-1)}- \xi_n(u,h))$,
\item[(iv)] there exists a constant $k\in \C$ such that $h$ satisfies 
\beq\label{btn}
 {\rm (BT)}_{u,k}
 \quad \bca h_x^{(n-1)} = \xi_n(u,h)-k,\\ h_t= \eta_{n,j}(u,h).\eca 
\eeq
\een

\bthm\label{ok1} (\cite{TWc})  Suppose $E$ is a frame of a solution $u=\sum_{i=1}^{n-1} u_i e_{in}$ of the $j$-th $\an1$-KdV flow \eqref{dsj} such that $E(x,t,\l)$ is holomorphic for $\l$ in an open subset $\calO$ of $\C$.  Let $k\in \co$, ${\bf c_0}= (c_1, \ldots, c_{n-1}, -1)^t$ a constant vector in $\C^n$, $(v_1, \ldots, v_n)^t:= E(\cdot, \cdot, k)^{-1}({\bf c_0})$, and $h:= -\frac{v_{n-1}}{v_n}$.  Then  we have the following:
\ben
\item $h$ is a solution of {\rm (BT)}$_{u,k}$, i.e., \eqref{btn}.
\item All solutions of \eqref{btn} are obtained this way.
\item $\ti u$ defined by \eqref{nm1} is a solution of \eqref{dsj}, (we will denote $\ti u$ as $h\ast u$).
\item $\det(f_{u,h})= (-1)^{n-1} (\l-k)$.
\item $\ti E= Ef_{u,h}^{-1}$ is a frame of $\ti u$ and $\ti E(x,t,\l)$ is holomorphic for $\l\in \calO\bh \{k\}$.
\item Let $C(\l)= e_{1n}(\l-k) + b +\sum_{i=1}^{n-1} c_i e_{i+1, n}$. Then 
$$\ti E(x,t,\l) = C(\l) E(x,t,\l) f_{u,h}(x,t,\l)^{-1}$$
is a frame of $\ti u$ that is  holomorphic for all $\l\in \calO$. 
\een
\ethm

As a consequence of Theorems  \ref{cfagdnc}, \ref{cfagdn}, and \ref{ok1}, we obtain BTs for the $j$-th central affine flow \eqref{ncfj}:

\bthm {\rm [BT for the $j$-th central affine curve flow with $k\not=0$]}\label{nx} \  

\ni
Let $\g(x, t)$ be a solution of the $j$-th  central affine curve flow \eqref{ncfj} on $\calM_n(\R)$, $g(\cdot, t)$ and $u(\cdot, t)$ the central affine moving frame and curvature of $\g(\cdot, t)$ respectively.  Let $k$ be a non-zero real constant, and $d(k)\in GL(n,\R)$ such that $\det(d(k))= (-1)^n k$. Let $h$ be a solution of {\rm (BT)}$_{u,k}$ (i.e., \eqref{btn}), and 
$$\ti g(x,t)=d(k) g(x,t) f_{u,h}(x,t,0)^{-1}.$$
Then
$h\ast \g(x,t)=\ti g(x,t) e_1$
is again a solution of \eqref{ncfj} and $\ti g(\cdot, t)$ is the central affine moving frame along $h\ast\g(\cdot, t)$. 
\ethm

\begin{proof} By Theorem \ref{cfagdn}, $u$ is a solution of \eqref{dsj}.  Since $h$ is a solution of (BT)$_{u,k}$, by Theorem \ref{ok1} (4) we have $\det(f_{u,h})= (-1)^{n-1} (\l-k)$ and $\hat g(x,t)= g(x,t)f_{u,h}^{-1}(x,t, 0)$ is a parallel frame of the Lax pair for $\ti u$ at $\l=0$.  So  $\det(\hat g)= \det(gf_{u,h}^{-1}(x,t,0))= (-1)^n k^{-1}$.  Since $d(k)$ is a constant matrix, $\ti g=d(k) \hat g$ is a frame of the Lax pair of $\ti u$ at $\l=0$ with determinant $1$.   Then this theorem follows from Theorems \ref{cfagdnc} and \ref{ok1}.
\end{proof}

Given an $n\times n$ matrix $M$, we use $M^\sharp$ to denote the {\it cofactor matrix of $M$\/}, i.e., the $ij$-th entry of $M^\sharp$ is equal to $(-1)^{i+j}$ times the  determinant of the $(n-1)\times (n-1)$ matrix obtained by crossing out the $j$-th row and the $i$-th column of $M$.  Then $MM^\sharp= \det(M)\I_n$.  

\bthm {\rm [BT for the $j$-th central affine curve flow with $k=0$]}  \label{pb} \ 

\ni Let $\g$, $u$, and $g$ be as in Theorem \ref{nx}, and $E(x, t, \l)$ the frame for the solution $u$ of \eqref{dsj} satisfying $E(0, 0, 0)=g(0, 0)$. Let $E_1(x,t)= \frac{\p}{\p\l}\big|_{\l=0} E(x,t,\l)$, ${\bf c}_0=(c_1, \ldots, c_{n-1}, -1)^t$ a constant in $\R^n$, $(v_1, v_2, \cdots, v_n)^t:= g^{-1}{\bf c}_0$, and $h=-\frac{v_{n-1}}{v_n}$. Then
$$ \ti{\g}=(e_{1n}g+AE_1)f_{u,h}^{\sharp}(x, t, 0)e_1$$ is a solution of \eqref{ncfj},  where $A= b+\sum_{i=1}^{n-1} c_i e_{i+1, n}$ and  $f_{u,h}^\sharp$ is the cofactor matrix of $f_{u,h}$.
\ethm

\begin{proof} Since both $E(x,t,0)$ and $g(x,t)$ are solutions of 
$$\bca g^{-1}g_x= b+u,\\ g^{-1}g_t= Z_j(u,0),\eca$$ with the same initial data at $(0,0)$, we have $E(x,t,0)= g(x,t)$. 

 Let $C(\l)= J+ \sum_{i=1}^{n-1} c_i e_{i+1, n}$. By Theorem \ref{ok1}, 
$$F(x,t,\l)= C(\l) E(x,t,\l) f_{u,h}^{-1}(x,t,\l)=\frac{(-1)^{n-1}}{\l}\, C(\l) E(x,t,\l) f_{u,h}^\sharp (x,t,\l)$$
is a frame for $h\ast u$ and is holomorphic at $\l=0$.  So
\beq\label{on}
F(x,t,0) =(-1)^{n-1} \frac{\p}{\p\l}\big|_{\l=0} C E f^\sharp_{u,h}.
\eeq
 By \eqref{on}, we get
 $$F(x,t, 0)=(e_{1n}g+AE_1)f_{u, h}^{\sharp}(x, t, 0)+Ag\left(\frac{\p }{\p \l}\big|_{\l=0} f_{u, h}^{\sharp}\right).$$
But the first column of $f_{u, h}^{\sharp}$ is independent of $\l$, hence the last term in the above formula is zero and
$$F(x,t, 0)=(e_{1n}g+AE_1)f_{u, h}^{\sharp}(x, t, 0).$$

Note that $\det(E(x,t,0))= \det(g(x,t))=1$. 
By Theorem \ref{ok1},  $\det(f_{u,h})= (-1)^{n-1}\l$. Note that  $\det(C(\l))= (-1)^{n-1} \l$. Hence $\det(F(x, t, 0))=1$. By Theorem \ref{cfagdnc}, $F(\cdot, \cdot, 0)e_1$ is a solution of the curve flow \eqref{ncfj}. 
\end{proof}
 
\bthm (\cite{TWc}) {\rm [Permutability for BT of the $\an1$-KdV hierarchy]}\  \label{pa}

\ni Let $u$ be a solution of \eqref{dsj}, $k_1, k_2\in \C$ constants,  $h_i$ solutions of {\rm (BT)}$_{u,k_i}$ \eqref{btn}, and $h_i\ast u$ the solution of \eqref{dsj} constructed from $u$ and $h_i$ for $i=1,2$ as in Theorem \ref{ok1}. Suppose $h_1 \neq h_2$. 
Set
$$
\bca
\ti{h}_1=h_1+\frac{(h_1-h_2)_x}{h_1-h_2}, \\
\ti{h}_2=h_2+\frac{(h_1-h_2)_x}{h_1-h_2}.
\eca
$$
Then
\ben
\item[(i)] $\ti h_1$ is a solution of {\rm (BT)}$_{h_2\ast u, k_1}$ and $\ti h_2$ is a solution of {\rm (BT)}$_{h_1\ast u, k_2}$,
\item[(ii)] $\ti h_1\ast (h_2 \ast u)= \ti h_2\ast(h_1\ast u)$. 
\een
\ethm

As a consequence of Theorems \ref{cfagdnc} and \ref{pa}, we obtain the following. 

\bthm {\rm [Permutability for the $j$-th central affine curve flow]} \par

Let $\g$, $u$, $g$, and $d(k)$ be as in Theorem \ref{nx}, $k_1, k_2, h_1, h_2, \ti h_1, \ti h_2,$  as in Theorem \ref{pa}, and 
$$g_i(x,t)= d(k_i) g(x,t) f_{u, h_i}(x,t,0)^{-1}, \quad \g_i= g_i e_1$$
for $i=1,2$, where $e_1= (1,0, \ldots, 0)^t$. 
Then 
$\ti h_2\ast \g_1= \ti h_2\ast \g_2$ is again a solution of \eqref{ncfj}, i.e.,  
$$ d(k_2) g_1(x,t) f_{u,\ti h_2}(x,t,0)^{-1}e_1= d(k_1) g_2(x,t) f_{u, \ti h_1} (x,t,0)^{-1}e_1$$
is again a solution of \eqref{ncfj}. 
\ethm

Next we write down formulas of BTs for the second central affine curve flow on $\calM_3(\R)$ in explicit forms.  
The second central affine curve flow on $\calM_3(\R)$ is 
\beq\label{cu} 
\g_t= -\frac{2}{3} u_2 \g + \g_{xx},
\eeq
and the second $A_2^{(1)}$-KdV flow is 
\beq\label{kdv2}
\begin{cases}
(u_1)_t = (u_1)_{xx}-\frac{2}{3}(u_2)_{xxx}+\frac{2}{3}u_2(u_2)_x,\\
(u_2)_t = -(u_2)_{xx}+2(u_1)_x.
\end{cases}
\eeq
System (BT)$_{u,k}$ for the second flow \eqref{kdv2} is 
\beq\label{bt3}
\bca 
h_{xx}= -u_1+(u_2)_x + hu_2 -3hh_x - h^3-k,\\
h_t=\frac{2}{3}(u_2)_x-h_{xx}-2hh_x.
\eca
\eeq
and the new solution  $h\ast u=\ti u_1 e_{13}+ \ti u_2 e_{23}$ is given by 
\beq\label{ct}
\bca \ti u_1= u_1- (u_2)_x + 3hh_x,\\ \ti u_2= u_2-3h_x.\eca
\eeq
Moreover, if $E$ is a frame for $u$, then $Ef_{u,h}^{-1}$ is a frame for $h\ast u$, where
$$
f_{u,h}(x,t,\l)= e_{13}\l + \bpm h& h_x& u_1-(u_2)_x+ h_{xx} + 3hh_x\\ 1 & h & u_2- h_x\\ 0& 1&h\epm.
$$

We choose $d(k)= -k^{1/3}\I_3$ in Theorem \ref{nx} for $n=3$. Then we obtain:

\bcor\label{nxa}  Let $\g$ be a solution of \eqref{cu} on $\calM_3(\R)$,  $u=\sum_{i=1}^2 u_ie_{i3}$ the central affine curvature of $\g(\cdot, t)$, $k$ a non-zero real constant, and $h$ a solution of \eqref{bt3}.  Then 
$$h\ast \g= k^{-2/3} ((h^2+h_x -u_2) \g - h \g_x+ \g_{xx})$$
is a solution of \eqref{cu} with central affine curvature $\ti u_1, \ti u_2$ given by \eqref{ct}. 
\ecor

\bcor \label{btw}
Let $\g$ and $u$ be as in Corollary \ref{nxa}, $g(\cdot, t)$ the central affine moving frame along $\g(\cdot, t)$, and $E(x, t, \l)$ the frame of $u$ with $E(0, 0, \l)=g(0, 0)$. Let $E_1(x,t)= \frac{\p}{\p\l}\n_{\l=0} E(x,t,\l)$, $c_1, c_2\in \R$, $A=e_{21}+ e_{32} +\sum_{i=1}^2 c_i e_{i+1,3}$,  $v=(v_1, v_2, v_3)^t=g^{-1}(c_1, c_2, -1)^t$, and  $h=-\frac{v_2}{v_3}$. Then
\ben
\item[(i)] $\ti{\g}=(e_{13}g+AE_1)(h^2+h_x-u_2, -h, 1)^t$ is a solution of \eqref{cu}, 
\item[(ii)] $\ti{E}(x,t,\l)=(\l e_{13}+A)E(x,t, \l)f_{u, h}^{-1}$
is a frame for $\ti{u}=h\ast u$.
\een
\ecor

Note that $[\p_x+ J, \p_t + J^j]=0$ implies that $u=0$ is a solution of the $j$-th $\an1$-KdV flow \eqref{dsj} and 
$E(x,t,\l)= \exp(xJ+ tJ^j)$ is a frame of $u=0$.  Hence 
\beq\label{tg}
\g(x,t)= \exp(bx+ b^jt) e_1
\eeq
 is a solution of the $j$-th central affine curve flow \eqref{ncfj} with zero central affine curvatures and $g(x,t)= \exp(bx+ b^jt)$ is the central affine moving frame of $\g$.  So we can apply Theorem \ref{nx} and \ref{pa} repeatedly to obtain an infinitely many families of explicit solutions of \eqref{ncfj}.

 \beg[Explicit 1-soliton solutions] \par

We apply BT to the solution $u=0$ of the second $A_2^{(1)}$-KdV flow \eqref{kdv2} with $k=8\mu^3\in \C$ (for detail cf. \cite{TWc}) to obtain $1$-soliton complex valued solutions
$$
\bca
u_1=9\mu^3\sech^2(\sqrt{3}\mu(x-2\mu it))(\sqrt{3}\tanh(\sqrt{3}\mu(x-2\mu it))+i), \\
u_2=-9\mu^2\sech^2(\sqrt{3}\mu(x-2\mu it))
\eca
$$
 for \eqref{kdv2}.  
 Apply BT to $u=0$ with $k= -c^3$ with constant $c\in \R \bh 0$ to get the following real solutions of \eqref{kdv2},
$$
\bca u_1= 9c^3 \sec^2 (\sqrt{3} c (x+ 2ct))(1+\sqrt{3}\tan (\sqrt{3} c (x+ 2ct))),\\
u_2= 9c^2 \sec^2 (\sqrt{3} c (x+ 2ct)).\eca
$$
We have seen that $\g(x,t)=(1, x, \frac{1}{2} x^2+t)^t$ is a solution of  \eqref{cu} with zero central affine curvatures. 
So we apply Theorem \ref{nx} and obtain the following solutions of the second central affine curve flow \eqref{cu}:
$$\ti \g=  \bpm 2(\xi-1)\\ 2x(\xi-1)+c^{-1}(\xi+1)\\ (x^2+2t)(\xi-1)+c^{-1}x(\xi+1)+ c^{-2}\epm,$$
where $c\in \R \bh 0$ is a constant and $\xi(x,t)= \sqrt{3}\tan(\sqrt{3}c (x+2ct))$. 
\eeg

\beg \label{rss}{\rm [Rational solutions for $n=3$]}\par

We apply BT with parameter $k=0$ to the vacuum solution $u=0$ of \eqref{kdv2}. Recall that $\g(x, t)=(1, x, \frac{x^2}{2}+t)^t$ is a solution of \eqref{cu} with zero central affine curvatures and $E(x, t, \l)=e^{xJ+tJ^2}$ is a frame of the trivial solution of \eqref{kdv2}. The constant term of $E(x,t,\l)$ as a power series in $\l$ is 
$$g(x,t)= E_0(x,t)=\exp(bx+ b^2 t)= \bpm 1&0&0\\ x& 1&0\\ \frac{x^2}{2} + t & x& 1\epm,$$ which is the central affine moving frame along $\g$. A direct computation implies that the coefficient of $\l$ of $E(x, t, \l)$ is 
$$E_1(x, t)=\bpm xt+\frac{1}{6}x^3 & t+\frac{1}{2}x^2 & x \\ \frac{1}{2}t^2+\frac{1}{2}x^2t+\frac{1}{24}x^4 & xt+\frac{1}{6}x^3 & t+\frac{1}{2}x^2 \\ \frac{1}{2}xt^2+\frac{1}{6}x^3t+\frac{1}{5!}x^5 & \frac{1}{2}t^2+\frac{1}{2}x^2t+\frac{1}{24}x^4 & xt+\frac{1}{6}x^3\epm.$$
We apply Theorem \ref{pb} to $u=0$ and $v_{0}=(a_1, a_2, 1)^t$ to see that
$$\ti{\g}=\bpm (h^2+h_{xx})(\frac{x^2}{2}+t) -hx+1 \\ (h^2+h_{xx})(\frac{1}{6}x^3+xt)-h(\frac{1}{2}x^2+t)+x \\ (h^2+h_{xx})(\frac{1}{2}t^2+\frac{1}{2}x^2t+\frac{1}{24}x^4)-h(\frac{1}{6}x^3+xt)+\frac{1}{2}x^2+t\epm$$
is a rational solution of \eqref{cu}, where 
$$h=\frac{a_1x-a_2}{1+a_1(\frac{x^2}{2}-t)-a_2x}.$$ The central affine curvatures for $\ti \g$ are
$$
\bca
u_1=-\frac{3(a_1x-a_2)(\frac{1}{2}a_1^2x^2-a_1a_2x+a_1^2t+a_2^2-a_1)}{(1+a_1(\frac{x^2}{2}-t)-a_2x)^3},\\
u_2=\frac{3(\frac{1}{2}a_1^2x^2-a_1a_2x+a_1^2t+a_2^2-a_1)}{(1+a_1(\frac{x^2}{2}-t)-a_2x)^2}.
\eca
$$

If we apply Theorem \ref{pb} or the permutability formula repeatedly, then we obtain infinitely many families of rational solutions for the second central affine curve flow \eqref{cu} on $\cm_3(\R)$.
\eeg

\bs

\section{Bi-Hamiltonian structure for the $\an1$-KdV hierarchy}\label{go}

We first review the bi-Hamiltonian structure and commuting conservation laws for the $\an1$-KdV hierarchy (cf. \cite{Dic03}, \cite{DS84}). Then we write down the formulas for the Poisson operators in terms of the operator $P_u$ defined in Definition \ref{gp} and compute the kernel (i.e., Casimirs) of the these operators.   

The group $C^\infty(S^1, N_n^+)$ acts on $C^\infty(S^1, \calB_n^+)$ by gauge transformation
$$g\ast q= g(b+q) g^{-1}- g_x g^{-1} -b,$$
or equivalently, 
$$g\ast (\p_x+b+q)= g(\p_x+b+q) g^{-1}= \p_x+ b+ g\ast q.$$  It was proved in \cite{DS84} (cf. also in \cite{TU14b}) that $C^\infty(S^1, V_n)$ is a cross section of the gauge action $C^\infty(S^1,N_n^+)$ on $C^\infty(S^1, \calB_n^+)$, i.e., given $q\in C^\infty(S^1, \calB_n^+)$ there exist unique $\D\in C^\infty(S^1, N_n^+)$ and $u\in C^\infty(S^1, V_n)$ such that $\D\ast q= u$.  In other words, each gauge orbit of $C^\infty(S^1, N_n^+)$ meets $C^\infty(S^1, V_n)$ exactly once. So $C^\infty(S^1, V_n)$ is isomorphic to the orbit space $\frac{C^\infty(S^1,\calB_n^+)}{C^\infty(S^1, N_n^+)}$. 

   Let $\li\, , \ri$ be the bi-linear form on $C^\infty(S^1, sl(n,\R))$ defined by
 \beq\label{qa}
 \li y_1, y_2\ri=\oint \tr(y_1(x) y_2(x))\rd x.
 \eeq The gradient of $\calF:C^\infty(S^1, \calB_n^+)\to \R$ at $q\in C^\infty(S^1, \calB_n^+)$ is the unique element $\K \cf(q)$ in $C^\infty(S^1, \calB_n^-)$ defined by
 $$\rd \calF_q(y)=\li \K \calF(q), y\ri$$
 for all $y\in C^\infty(S^1, \calB_n^+)$.  
  
 The two Poisson structures on $C^\infty(S^1, \calB_n^+)$ given in \cite{DS84} are:
 \begin{align}
 \{\calF_1, \calF_2\}_1(u)&= \li [e_{1n}, \K \calF_1(u)], \K \calF_2(u)\ri, \label{sx1}\\
 \{\calF_1, \calF_2\}_2(u)&= \li [\p_x+ b+ u, \K\calF_1(u)], \K \calF_2(u),\label{sx2}\ri.
 \end{align}
 These two Poisson structures are invariant under the action of $C^\infty(S^1, N_n^+)$.  So they induce two Poisson structures on $C^\infty(S^1, V_n)$. 

The gradient $\K F(u)$ of a functional $F$ on $C^\infty(S^1, V_n)$ at $u$  is the unique element in $C^\infty(S^1, V_n^t)$ satisfying
 $$\rd F_u(v)=\li \K F(u), v\ri$$  for all $v\in C^\infty(S^1, V_n)$. 
 
 Recall that $C^\infty(S^1, V_n)$ is isomorphic to the orbit space $\frac{C^\infty(S^1, \calB_n^+)}{C^\infty(S^1, N_n^+)}$. So given a functional $F$ on $C^\infty(S^1, V_n)$, there is a unique $C^\infty(S^1, N_n^+)$-invariant function 
  $\ti F$ on $C^\infty(S^1, \calB_n^+)\to \R$ whose restriction to $C^\infty(S^1, V_n)$ is $F$, i.e., $\ti F(q)= F(u)$ if $u\in C^\infty(S^1, V_n)$ lies in the same $C^\infty(S^1, N_n^+)$-orbit as $q$.  
 The following Proposition gives the relation between $\K \ti F(u)$ and $\K F(u)$ for $u\in C^\infty(S^1, V_n)$. 
 
\bprop Let $F$ be a functional on $C^\infty(S^1, V_n)$, and $\ti F$ the functional on $C^\infty(S^1,\calB_n^+)$ invariant under  $C^\infty(S^1, N_n^+)$ defined by $F$.  Then
$$\K \ti F(u)= P_u(\K F(u)),$$
where $u\in C^\infty(S^1, V_n)$ and $P_u$ is the operator defined Definition \ref{gp}.  
\eprop

\begin{proof} Note that the infinitesimal vector field $\ti \xi$ defined by $\xi$ in $ C^\infty(S^1, \calN_n^-)$ for the gauge action is 
$$\ti\xi(q)= -[\p_x+b+q, \xi].$$
where $q\in C^\infty(S^1, \calB_n^+)$. By assumption, $\ti F(f\ast q)= \ti F(q)$ for all $q\in C^\infty(S^1, \calB_n^+)$ and $f\in C^\infty(S^1, N_n^+)$. So
$\rd \ti F_q(\ti\xi(q))=0$. But 
$$ \rd \ti F_q(\ti\xi(q))=\li \K \ti F(q), \ti\xi(q)\ri= -\li \K \ti F(q), [\p_x+b+q, \xi]\ri = \li [\p_x+b+q, \K \ti F(q)], \xi\ri$$
for all $\xi\in C^\infty(S^1, \calN_n^+)$. So 
$$
[\p_x+b+u, \K \ti F(q)]\in C^\infty(S^1, \calB_n^-).
$$
  
To prove $\K \ti F(u)= P_u(\K F(u))$ for $u\in C^\infty(S^1, V_n)$ is equivalent to prove
\beq\label{nd}
\rd\ti F_u(y)= \li P_u(\K F(u)), y\ri
\eeq
for all $y\in C^\infty(S^1, \calB_n^+)$. 
  
We first prove \eqref{nd} for $y\in C^\infty(S^1, V_n)$. Given $u, v\in C^\infty(S^1,V_n)$, we have  
$$\rd F_u(v)= \li \K F(u), v\ri = \rd \ti F_u(v) =\li \K \ti F(u), v\ri.$$
So $\li \K F(u)-\K \ti F(u), v\ri=0$ for all $v\in C^\infty(S^1, V_n)$.  This implies that 
$$\pi_0(\K \ti F(u))= \K F(u),$$
where $\pi_0:sl(n,\R)\to V_n^t$ is the canonical projection. By definition of $P_u$, we have $\pi_0(P_u(\K F(u)))= \K F(u)$. So we obtain $\rd \ti F_u(v)= \rd F_u(v)= \li P_u(\K F(u)), v\ri$, i.e., \eqref{nd} is true for $y\in C^\infty(S^1, V_n)$.  
   
Since $C^\infty(S^1, V_n)$ is a cross section of the gauge action of $C^\infty(S^1, N_n^+)$ on $C^\infty(S^1, \calB_n^+)$, the tangent space of $C^\infty(S^1, \calB_n^+)$ at $u\in C^\infty(S^1, V_n)$ can be written as a direct sum of $C^\infty(S^1, V_n)$ and the tangent space of the $C^\infty(S^1, N_n^+)$-orbit at $u$.  Since $\ti F$ is invariant under $C^\infty(S^1, N_n^+)$, we have $\rd \ti F_u(\ti\xi(u))=0$ for all $\xi\in C^\infty(S^1, \calN_n^+)$.  So 
\begin{align*}
\li P_u(\K F(u)), \ti\xi(u)\ri &=\li P_u(\K F(u)), -[\p_x+b+u, \xi]\ri\\ 
&= \li [\p_x+b+u, P_u(\K F(u))], \xi\ri.
\end{align*}
By definition of $P_u$, we have $[\p_x+b+u, P_u(\K F(u))]\in C^\infty(S^1, V_n)$.  Since $\xi\in C^\infty(S^1, \calN_n^+)$, we conclude that 
$\li P_u(\K F(u)),\ti \xi(u)\ri=0$. This proves $\rd \ti F_u(\ti \xi(u))=\li P_u(\K F(u)), \ti\xi(u)\ri=0$.  So \eqref{nd} is true for $y$ in the tangent space of $C^\infty(S^1, N_n^+)$-orbit at $u$.  This completes the proof.    
\end{proof}

So the Poisson structures on $C^\infty(S^1, V_n)$ induced from \eqref{sx1} and \eqref{sx2} are given as follows: 
\begin{align*}
& \{F_1, F_2\}_1(u)= \li [e_{1n}, P_u(\K F_1(u))], P_u(\K F_2(u))\ri, \\
&\{F_1, F_2\}_2(u)= \li [\p_x+b+u, P_u(\K F_1(u))], P_u(\K F_2(u))\ri, 
\end{align*} 
where $P_u:C^\infty(S^1,V_n^t)\to C^\infty(S^1, sl(n,\R))$ is defined in Definition \ref{gp}. 

Let  $$(J_i)_u:C^\infty(S^1,V_n^t)\to C^\infty(S^1, V_n)$$ be the Poisson operator corresponding to $\{\, , \}_i$ at $u$ for $i=1, 2$, i.e., $(J_i)_u$ is defined by
$$
\{F_1, F_2\}_i(u)=\li (J_i)_u(\K F_1(u)), \,\K F_2(u)\ri.
$$
Then the Hamiltonian equation for a functional $H:C^\infty(S^1, V_n)\to \R$ with respect to $\{\, , \}_i$ is 
$$u_t=(J_i)_u(\K H(u)).$$

Next we compute the formula for the Poisson operator $J_1$. 

\bprop\label{hr} The Poisson operator $(J_1)_u:C^\infty(S^1, V_n^t)\to C^\infty(S^1, V_n)$  is of the form
$(J_1)_u(\xi)= \sum_{i=1}^{n-1} (L_i)_u(\xi) e_{in}$ with 
$$
(L_i)_u(\xi)= n(\xi_{n-i})_x + k_i(\xi_1, \ldots, \xi_{n-i-1}),
$$
where $\xi= \sum_{i=1}^{n-1} \xi_i e_{ni}$ and $k_i$'s are linear differential operators with differential polynomials of $u$ as coefficients. 
\eprop

\begin{proof} By Theorem \ref{jw}, entries of $P_u(v)$ are differential polynomials of $u$ and $v$. So we can use integration by parts to compute $(J_1)_u$. We proceed as follows: 
Let $u=\sum_{i=1}^{n-1} u_i e_{in}$, 
and 
\begin{align*}
&\xi=\sum_{i=1}^{n-1} \xi_i e_{ni}:= \K F_1(u), \quad \eta=\sum_{i=1}^n \eta_i e_{ni}:= \K F_2(u), \\
& C=(C_{ij})= P_u(\xi), \quad D= (D_{ij})= P_u(\eta).
\end{align*} 
Then we have 
$$\{F_1, F_2\}_1(u)= \li [e_{1n}, C], D\ri = \oint \sum_{i=1}^n C_{ni} D_{i1} - C_{i1} D_{ni}  \rd x.$$

Let $g:\R\to SL(n,\R)$ be a solution of $g^{-1}g_x= b+u$, $\g$ the first column of $g$. Then $\g\in \calM_n(\R)$, $g$ is the central affine moving frame along $\g$, and $u=\Psi(\g)$,  By Corollary \ref{sy}, there is a $\d\g$ tangent to $\calM_n(\R)$ at $\g$ such that $C= g^{-1}\d g$, where $\d g=(\d\g, \ldots, (\d\g)_x^{(n-1)})$. Since $\d\g=\sum_{i=1}^n C_{i1} \g_x^{(i-1)}$ is tangent to $\calM_n(\R)$, it follows from Theorem  \ref{jw} that there exists a differential polynomial $\phi_0$ such that
$$C_{11}= \phi_0(C_{21}, \ldots, C_{n1})= f(\xi_1, \ldots, \xi_{n-1}).$$

By definition, $C_{ni}= \xi_i$ and $D_{ni}= \eta_i$. By Theorem  \ref{jw}, there is a differential polynomial $f_n$ such that
$$C_{nn}=f_n(\xi_1, \ldots, \xi_{n-1}), \quad D_{nn}= f_n(\eta_1, \ldots, \eta_{n-1}).$$
We then use integration by parts to write down the Poisson operator $(J_1)_u$.
To get $(L_j)_u(\xi)$, we only need to calculate the terms involving $\eta_j= D_{nj}$ in $\sum_{i=1}^n C_{i1}D_{ni} - C_{ni} D_{i1}.$
We use \eqref{he} to compute these terms as follows:
\begin{align*}
& \quad D_{nj}C_{j1}+C_{n1}(D_{nn}-D_{11})-\sum_{i=1}^{n+1-j}C_{ni}D_{i1} \\
& =\eta_jC_{j1}+\xi_1\left(\sum_{i=1}^{n-1}D_{i+1, i}'\right)-C_{n, n+1-j}D_{n+1-j, 1}-\sum_{i=1}^{n-j}C_{ni}D_{i1} \\
&=\eta_j(\xi_{n+1-j}-(n-j)\xi_{n-j}'+\phi_{n-j}(u, \xi_1, \cdots, \xi_{n-j-1}))+\xi_1\sum_{i=1}^{n-1}D_{i+1, i}' \\
& \quad -\xi_{n+1-j}(\eta_{j}-(j-1)\eta_{j-1}'+\phi_{j-1}(u, \eta_1, \cdots, \eta_{j-2}))-\sum_{i=1}^{n-j}\xi_{i}D_{i1}.
\end{align*}
Note that $\xi_1\sum_{i=1}^{n-1}D_{i+1, i}'$ only depends on $\xi_1$, $\eta_1, \cdots, \eta_{n-1}$,  and $\xi_iD_{i1}$ is a differential polynomial in $u$, $\xi_i$ and $\eta_1, \cdots, \eta_{n+1-i}$ for each $i$. Therefore, to consider the term $\xi_{n-j}$ in the the coefficient of $\eta_j$ in $\sum_{i=1}^{n-j}\xi_{i}D_{i1}$, we only need to calculate $\xi_{n-j}D_{n-j, 1}$. Again, use $D_{n-j, 1}=\eta_{j+1}-j\eta'_{j}+\phi_{j}(u, \eta_1, \cdots, \eta_{j-1})$ and 
 integration by parts to see that the coefficients of $\eta_j$ is 
$-n\xi_{n-j}'$ plus a differential operator depending on $u, \xi_1. \cdots, \xi_{n-j-1}$.
\end{proof}

\bcor\label{hr1}
The dimension of the kernel of $(J_1)_u$ is $n-1$. 
\ecor

\begin{proof}
If $(J_1)_u(\xi)=0$ with $\xi=\sum_{i=1}^{n-1}\xi_ie_{ni}$, then Proposition \ref{hr} implies that $(L_{n-1})_u(\xi)=n(\xi_1)_x=0$. So $\xi_1=c_1$ a constant. Use $(L_{n-2})_u(\xi)= n(\xi_2)_x+ k_{n-2}(\xi_1)=0$ to see that $\xi_2= -\frac{1}{n} k_{n-2}(c_1) x + c_2$ for some constant $c_2$. The corollary follows from induction. 
\end{proof}

Next we compute the formula for the second Poisson operator $J_2$. 

\bprop
The Poisson operator $(J_2)_u:C^\infty(S^1,V_n^t)\to C^\infty(S^1,V_n)$ is 
\beq\label{hn}
(J_2)_u(v)=[\p_x+b+u, P_u(v)].
\eeq
\eprop

\begin{proof}
By definition of $P_u$, we have $[\p_x+b+u, P_u(\K F_1(u))]\in C^\infty(S^1, V_n)$. So we have 
$$\li [\p_x+b+u, P_u(\K F_1(u))], P_u(\K F_2(u))\ri=\li [\p_x+b+u, P_u(\K F_1(u))], \K F_2(u)\ri.$$ Hence the second Poisson structure can be written as 
$$
\{F_1, F_2\}_2(u)=  \li [\p_x+b+u, P_u(\K F_1(u))], \K F_2(u)\ri,
$$
which proves that $(J_2)_u(v)= [\p_x+b+u, P_u(v)]$ is the second Poisson operator.
\end{proof}

In the following examples, we compute explicit formulas for $J_1$ and $ J_2$ for small $n$.

\beg For $n=2$, write $u=qe_{21}$, $\K F_1(u)=\xi e_{12}$, and $\K F_2(u)=\eta e_{21}$, then 
\begin{align*}
& P_u(\K F_1(u))=\bpm -\frac{1}{2}\xi_x & -\frac{1}{2}\xi_{xx}+q\xi  \\
\xi & \frac{1}{2}\xi_x\epm, \\
& P_u(\K F_2(u))=\bpm -\frac{1}{2}\eta_x & -\frac{1}{2}\eta_{xx}+q\eta  \\
\eta & \frac{1}{2}\eta_x\epm.
\end{align*}
So we get
\begin{align*}
&\{F_1, F_2\}_1(u)=2\oint \xi'\eta \rd x,  \\ 
&\{F_1, F_2\}_2(u)=-\oint (\frac{1}{2}\xi_x^{(3)}-2q\xi_x-q_x\xi)\eta \rd x, 
\end{align*} and the corresponding Poisson operators are
\begin{align*}
& (J_1)_u(\xi e_{21})=2\xi_xe_{12}, \\
&(J_2)_u(\xi e_{21})=(-\frac{1}{2}\xi_x^{(3)}+2q\xi_x+q_x\xi) e_{12}.
\end{align*}
These are the known Poisson structures for the KdV hierarchy.
\eeg

\beg For $n=3$, write $u=u_1e_{13}+u_2e_{23}$, $\xi= \K F_1(u)=\xi_1 e_{31} + \xi_2 e_{32}$, $\eta= \K F_2(u)= \eta_1 e_{31} + \eta_2 e_{32}$, $C=P_u(\xi)=(C_{ij})$, and $D= P_u(\eta)=(D_{ij})$. Use the algorithm given in the proof of Theorem \ref{jw}  to compute $P_u(v)$ for $v= v_1 e_{31}+ v_2 e_{32}$ and  get
\beq\label{hf}
P_u(v)=\bpm 
-v_2'+\frac{2}{3} v_1'' -\frac{2}{3} u_2 v_1 & p_{12} &p_{13}\\ v_2-v_1' & p_{22} &p_{33}\\ v_1 & v_2 & v_2'-\frac{1}{3} v_1'' +\frac{1}{3} u_2 v_1
\epm,\eeq
where 
\begin{align*}
p_{12}&=-v_2''+\frac{2}{3}v_1^{(3)}-\frac{2}{3}(u_2v_1)'+u_1v_1, \\
p_{13}&=-v_2^{(3)}+\frac{2}{3}v_1^{(4)}-\frac{2}{3}(u_2v_1)''+(u_1v_1)'+u_1v_2,\\
p_{22}&=-\frac{1}{3}v_1''+\frac{1}{3}u_2v_1, \\
p_{23}&=-v_2''+\frac{1}{3}v_1^{(3)}-\frac{1}{3}(u_2v_1)'+u_2v_2+u_1v_1.
\end{align*}
Here we use $y'$ to denote $y_x$.
Then integration by part gives
$$\{F_1, F_2\}_1(u)= \oint \sum_{i=1}^3 C_{3i} D_{i1} - C_{i1}D_{3i}\rd x = 3\oint (\xi_1' \eta_2 + \xi_2' \eta_1)\rd x.$$
Hence 
$$(J_1)_u(\xi_1e_{31} + \xi_2 e_{32})= 3(\xi_2' e_{13} + \xi_1'e_{23}).$$ 
This formula was also obtained in \cite{CIM13}. 

We use \eqref{hn} and a direct computation to see that
\begin{align*}
(J_2)_u(\xi)&=(C_{13}'+u_1(C_{33}-C_{11})-u_2C_{12})e_{13} \\ & + (C_{23}'+C_{13}+u_2(C_{33}-C_{22})-u_1C_{21})e_{23}.
\end{align*}
This gives 
$$(J_2)_u(\xi_1 e_{31}+ \xi_2e_{32})= (A_1)_u(\xi) e_{13}+ (A_2)_u(\xi) e_{23},$$
where 
\begin{align*}
(A_1)_u(\xi)&=\frac{2}{3}\xi_1^{(5)}-\xi_2^{(4)}-\frac{2}{3}(u_2\xi_1)^{(3)}-\frac{2}{3}u_2\xi_1^{(3)}+u_1''\xi_1+2u_1'\xi_1' \\& \quad+3u_1\xi_2'+u_1'\xi_2 +u_2\xi_2''+\frac{2}{3}u_2(u_2\xi_1)'.\\
(A_2)_u(\xi) &=\xi_1^{(4)}-2\xi_{2}^{(3)}-(u_2\xi_1)''+2u_2\xi_2'+u_2'\xi_2+u_1\xi_1'+2u_1'\xi_1.
\end{align*}
\eeg

\beg
For $n=4$, let $u=\sum_{i=1}^3u_ie_{i4}$, $\xi=\K F_1(u)=\sum_{i=1}^3\xi_ie_{4i}$, $\eta=\K F_2(u)=\sum_{i=1}^3\eta_ie_{4i}$, $P_u(\xi)=(C_{ij})$ and $P_u(\eta)=(D_{ij})$. Then the algorithm given in the proof of Theorem \ref{jw} gives
$$P_u(v)=\bpm a_{11} & \ast & \ast & \ast \\ v_3-2v_2'+v_1''-u_3v_1 & \ast & \ast & \ast \\ v_2-v_1' & \ast & \ast & \ast \\ v_1 & v_2 & v_3 & a_{44}
\epm,$$
where $v=\sum_{i=1}^3 v_i e_{4i}$, and 
\begin{align*}
a_{11}&=-\frac{1}{4}(3v_1^{(3)}-8v_2''+6v_3'-3(u_3v_1)'+3u_2v_1+2u_3v_2), \\
a_{44}&=\frac{1}{4}v_1^{(3)}-v_2''+\frac{3}{2}v_3'-\frac{1}{4}(u_3v_1)'+\frac{1}{4}u_2v_1+\frac{1}{2}u_3v_2.
\end{align*}
 So the first Poisson structure is 
\begin{align*}
&\{F_1, F_2\}_1(u)=\oint \sum_{j=1}^{4}(C_{4j}D_{j1} - C_{j1}D_{4j})\rd x \\
&=4\oint (\xi_3'-\frac{1}{2}\xi_2''+\frac{1}{2}\xi_1^{(3)}-\frac{1}{4}(2u_3\xi_1'+u_3'\xi_1))\eta_1+(\xi_2'+\frac{1}{2}\xi_1'')\eta_2+\xi_1'\eta_3 \rd x.
\end{align*}
Therefore, 
$$
(J_1)_u(\xi)=\bpm 4\xi_3'-2\xi_2''+2\xi_1'''-2u_3\xi_1'-u_3'\xi_1 \\ 4\xi_2'+2\xi_1'' \\ 4\xi_1'\epm.
$$
The formula for $(J_2)_u$ can be computed in a similar way as in the case $n=3$, but it is very long and we do  not include here. 
\eeg

Next we review the construction of conservation laws of the $\an1$-KdV hierarchy in \cite{DS84}.  

\bthm\label{ta} (\cite{DS84}) Given $u \in C^{\infty}(\R,V_n)$, then there exists a unique 
$T=\sum_{i=0}^{\infty}T_i\l^{-i}$
such that $T_0 \in C^{\infty}(\R, N_n^+)$, the first column of $T$ is
$e_1$, and
\beq\label{th}
T(\p_x+ J+ u)T^{-1}=\p_x+J+\sum_{i=0}^{\infty}f_i(u)J^{-i}, \ \ \
f_i(u) \in C^\infty(S^1, \R).
\eeq
Moreover, let $H_i$ the functional on $C^\infty(S^1,V_n)$ defined by 
\beq\label{sz}
H_i(u)=n\oint f_i(u)\rd x.
\eeq
Then we have
\ben
\item[(a)]
$\K H_j(u)= \pi_0(Y_{j,0}(u))= \pi_0(Z_{j,0}(u))$,
\item[(b)] the $j$-th $\an 1$-KdV flow \eqref{mgd0} is
 the Hamiltonian equation  of $H_j$ with respect to $\{ , \}_2$,
 \item[(c)] the $j$-th $\an 1$-KdV flow \eqref{mgd0} is
 is  the Hamiltonian
equation of $H_{n+j}$ with respect to $\{ , \}_1$, 
\een
where $\pi_0:sl(n,\R)\to V_n^t$ is the canonical projection and $Y_{j,0}(u)$ and $Z_{j.0}(u)$ are defined by \eqref{fm} and \eqref{by} respectively.
\ethm

\beg  We use \eqref{th} to compute $f_i$ explicitly.  
For example, for $n=3$ and $u=u_1 e_{13} + u_2 e_{23}$, equation \eqref{th} implies that
\begin{align*}
&H_1(u)= \oint u_2 \rd x,\quad   H_2(u)= \oint u_1 \rd x,\quad H_4(u)= -\frac{1}{3} \oint u_1u_2 \rd x,\\
&\K H_1(u)= e_{32}, \quad \K H_2(u)= e_{31}, \quad \K H_4(u)= -\frac{1}{3} u_2 e_{31}-\frac{1}{3} u_1 e_{32}.
\end{align*}  
Since $P_u(\K H(u))=g^{-1}\d g$ with $\pi_0(g^{-1}\d g)=\K H(u)$, where $\pi_0$ is the canonical projection defined in Definition \ref{gp}, use \eqref{hf} to see that 
\begin{align*}
&P_u(\K H_1(u))=\left(\begin{array}{ccc}0&0&u_1\\1&0&u_2
\\ 0&1&0\end{array}\right),\\
&Z_{2,0}(u)=P_u(\K H_2(u))=\left(\begin{array}{ccc}-\frac{2}{3}u_2&u_1-\frac{2}{3}u_2'&u_1'-\frac{2}{3} u_2''\\
0&\frac{1}{3}u_2&u_1-\frac{1}{3}u_2'\\1&0&\frac{1}{3}u_2\end{array}\right),\\
&Z_{4,0}(u)= P_u(\K H_4(u))=\left(\begin{array}{ccc}\frac{-2u_2^{(2)}+3u_1'+2u_2^2}{9}&\ast &\ast\\
\frac{u_2'-u_1}{3}&\frac{u_2''-u_2^2}{9}&0\\-\frac{u_2}{3}&-\frac{u_1}{3}&\frac{u_2''-3u_1'-u_2^2}{9}\end{array}\right).
\end{align*}

For general $n$, we have
\begin{align*}
H_1(u)&= \oint u_{n-1} \rd x,\\
 H_2(u)&=\oint u_{n-2} \rd x,\\
H_3(u) &= \oint u_{n-3}+\frac{n-3}{2n}u_{n-1}^2 \rd x.
\end{align*}
\eeg

  Next we use a similar computation as in \cite{Ter97} for the $n\times n$ AKNS hierarchy  to give another method to construct conservation laws for the $\an1$-KdV hierarchy from the solution $Y(u,\l)$ of \eqref{ff}. 

\bthm\label{rs} Given $u \in \C^{\infty}(S^1\times \R, V_n)$, let $Y(u,\l)$  be the solution of \eqref{ff}, $N(u,\l)= \l^{-1} Y^j(u,\l)$, and $\li \, , \ri$ the bi-linear form on $C^\infty(S^1,sl(n,\R))$ defined by \eqref{qa}. Then we have
$$
\left\li \frac{\p N(u, \l)}{\p \l}, \d u\right\ri= \d \li N(u, \l), e_{1n} \ri=  \d \li \l^{-1}Y^j(u, \l),  e_{1n} \ri
$$
for any variation $\d u$. 
\ethm

\begin{proof} Choose $M$ such that $M(\p_x+J+u)M^{-1}=\p_x+J$, i.e.,  
$$M^{-1}M_x+M^{-1}JM=J+u.$$
Then we have
$Y(u, \l)=M^{-1}JM$ and $Y^j(u, \l)=M^{-1}J^jM$. 

Set $\xi =M^{-1}M_\l$ and $\eta =M^{-1}\d M$. Direct computations give the following formulas:
\begin{align*}
& \d Y^j(u, \l) =[Y^j(u, \l), \eta], \\
& (Y^j)_x(u, \l)=[Y^j(u, \l), J+u], \\
& \frac{ \p Y^j(u, \l)}{\p \l }=[Y^j(u, \l), \xi]+M^{-1} (b^t)^{n-j}M,   \\
& \d u = \eta_x+[J+u, \eta], \\
& \xi_x+[J+u, \xi]=e_{1n}-M^{-1}e_{1n} M.
\end{align*}
Use the above formulas to compute to get
\begin{align*}
  \li \frac{\p Y^j}{\p \l}, \d u \ri  & = \li [Y^j, \xi]+M^{-1}(J^j)_\l M, \eta_x+[J+u, \eta] \ri \\
& = - \li [Y^j_x, \xi], \eta \ri -\li[Y^j(u, \l), \xi_x], \eta \ri \\
& \quad - \li (M^{-1}(J^j)_\l M)_x, \eta \ri+ \li M^{-1}(J^j)_\l M, [J+u, \eta] \ri \\
& = \li [Y^j, \eta], e_{1n} \ri+ \li \eta, M^{-1}([(J^j)_\l, J]-[e_{1n}, J^j])M \ri.
\end{align*}
Since $[J^j, J]=0$, we have $[(J^j)_\l, J]+ [J^j, e_{1n}]=0$. So we have
$$ \li \frac{\p Y^j(u, \l)}{\p \l}, \d u \ri = \d\li Y^j(u, \l), e_{1n} \ri.$$
Therefore 
\begin{align*}
\li \frac{\p N(u, \l)}{\p \l}, \d u\ri & = -\l^{-2}\li Y^j, \d u \ri+\l^{-1}\li \frac{\p Y^j}{\p \l}, \d u \ri \\
& = -\l^{-2}\li Y^j, \eta_x+[J+u, \eta] \ri+\l^{-1}\li \frac{\p Y^{j}}{\p \l}, \d u \ri  \\
& = \l^{-2}(\li (Y^j)_x, \eta \ri - \li Y^j, [J+u, \eta]\ri)+ \l^{-1}\li \frac{\p Y^j}{\p \l}, \d u \ri \\
&=\l^{-1}\li \frac{\p Y^j}{\p \l}, \d u \ri
\end{align*}
This proves the theorem.
\end{proof}

Since $Y^n(u,\l)= M^{-1}J^nM= \l\I_n$, we have $Y^{nk+j}(u,\l)= \l^k Y^j(u,\l)$.  So
$$Y_{nk+j,0}(u)= Y_{j, -k}(u).$$
Recall that  if the gradient of a functional at $u$ is $\pi_0(Y_{j,0}(u))$ then it is the Hamiltonian of the $j$-th $\an1$-KdV flow, $u_t=[\p_x+ b+u, P_u(\pi_0(Y_{j, 0}(u)))]$.   Therefore we have the following.

\bcor  Let $F_{j,i}:C^\infty(S^1, V_n)\to \R$ be the functional defined by
$$F_{j,i}(u)= -\frac{1}{i+1} \oint \tr(Y_{j, -(i+1)}(u) e_{1n}) \rd x,$$
where $Y_{j, -(i+1)}(u)$ is the coefficient of $\l^{-(i+1)}$ of $Y^j(u,\l)$ and $Y(u,\l)$ is as in Theorem \ref{rs}. Then $\K F_{i,j}(u)= \pi_0(Y_{j, -i}(u))=\pi_0(Y_{ni+j, 0}(u))$ and the Hamiltonian equation for $F_{j,i}$ with respect to $\{\, ,\}_2$ is the $(ni+j)$-th $\an1$-KdV flow.
\ecor

Next we compute the Casimirs of $\{\, ,\}_1$. 

\bthm\label{jr}
Let $H_j$ be the functional defined by \eqref{sz} for the $\an 1$-KdV hierarchy. Then $(J_1)_u(\K H_j(u))=0, 1 \leq j \leq n-1$. In other words, $H_1, \ldots, H_{n-1}$ are Casimirs of $\{\, ,\}_1$. Moreover, $\Ker((J_1)_u)$ is equal to the span of $\{\K H_1(u), \ldots, \K H_{n-1}(u)\}$. 
\ethm

\begin{proof}
Let $$Y(u,\l)= e_{1n}\l + Y_{1,0}(u) + Y_{1, 1}(u)\l^{-1} + \cdots$$ 
be the solution of \eqref{ff} for $u$, and $Y_{j,0}(u)$  the constant term of $Y(u,\l)^j$ as a power series in $\l$. 
We claim that $(J_1)_u(\pi_0(Y_{j,0}(u)))=0$ for $1\leq j\leq n-1$.
It follows from the expansion of $Y(u,\l)$  that we have
$$Y(u,\l)^j= (b^t)^{n-j} \l + Y_{j,0}(u) + Y_{j, 1}(u)\l^{-1} + \cdots, \quad 1\leq j\leq n-1.$$
Since $[\p_x+J+u, Y(u, \l)^j]=0$, $[\p_x+b+u, (b^t)^{n-j}]= [Y_{j,0}(u), e_{1n}]$.  Recall that $Z_{j,0}(u)= Y_{j,0}(u) -\zeta_j(u)$ for some unique $\zeta_j(u)\in C^\infty(S^1, \calN_n^+)$.  But $[\eta_j(u), e_{1n}]=0$. So we get 
$$[Z_{j,0}(u), e_{1n}]= [Y_{j,0}(u), e_{1n}]= [\p_x+ b+u, (b^t)^{n-j}].$$
By Corollary \ref{gpc}, $Z_{j,0}(u)= P_u(\pi_0(Z_{j,0}(u)))$. Let $\xi_1= \pi_0(Z_{j,0}(u))$, $\xi_2 \in C^{\infty}(S^1, V_n^t)$. Then we have
\begin{align*}& \li [Z_{j,0}(u), e_{1n}], P_u(\xi_2)\ri = \li [\p_x+b+u, (b^t)^{n-j}], P_u(\xi_2)\ri \\
&\quad = -\li (b^t)^{n-j}, [\p_x + b+u, P_u(\xi_2)]\ri,
\end{align*}
which is zero because $[\p_x + b+u, P_u(\xi_2)] \in V_n$ by definition of $P_u$ and $(b^t)^{n-j}$ is in $\calN_n^+$. This proves that $\li (J_1)_u(\pi_0(Z_{j,0}(u))), \xi_2\ri =0$ for all $\xi_2 \in V_n^t$. Hence $\pi_0(Z_{j,0}(u))$ lies in the kernel of $(J_1)_u$.  
By Theorem \ref{ta}, $\K H_j(u)=\pi_0(Z_{j, 0}(u))$ for $1 \leq j \leq n-1$. So we get $(J_1)_u(\K H_j(u))=0, 1 \leq j \leq n-1$.
\end{proof}

Below we derive some properties of the central affine curvature map $\Psi$  and use them to compute the kernel of $(J_2)_u$.

\bprop\label{fp} Let $\Psi:\calM_n(\R)\to C^\infty(\R, V_n)$ be the central affine curvature map and $u=\Psi(\g)= \sum_{i=1}^{n-1} u_i e_{in}$. Then 
$$
\rd \Psi_\g(\d \g)= [\p_x+ b+u, g^{-1}\d g]= (J_2)_u(\pi_0(g^{-1}\d g)), 
$$
where $g=(\g, \g_x, \ldots, \g_x^{(n-1)})$ is  the central affine frame for $\g$, $\pi_0:sl(n,\R)\to V_n^t$ is the canonical projection defined in Definition \ref{gp}, and $$\d g= (\d\g, (\d\g)_x, \ldots, (\d\g)_x^{(n-1)}).$$
\eprop

\begin{proof}
 Take variation of $g^{-1}g_x=b+u$ to get 
$$\d u= -(g^{-1}\d g) (b+u) + g^{-1}(\d g)_x.$$
Set $\eta= g^{-1}\d g$ and compute directly $\eta_x$ to get 
$\eta_x=-[b+u, \eta]+\d u$.  But $\Psi(\g)=g^{-1}g_x$, where $g= (\g, \g_x, \ldots, \g_x^{(n-1)})$. Hence $\rd \Psi_\g(\d\g)=\d u=  [\p_x+b+u, g^{-1}\d g]$.  
\end{proof} 

Recall that  $\Psi(\g_1)= \Psi(\g_2)$ if and only if there exists $c\in SL(n,\R)$ such that $\g_2= c\g_1$, where $\Psi$ is the central affine curvature map.   So we have

\bprop \label{fr}
Let $\Psi:\calM_n(S^1)\to C^\infty(S^1, V_n)$ be the central affine curvature map,  $\Psi(\g)=u$, and $g$ the central affine moving frame along $\g\in \calM_n(S^1)$. Then
\ben
\item $\Ker(\rd \Psi_\g)=\{c_0\g \,|\, c_0\in sl(n,\C)\}$,
\item $\Psi^{-1}(\Psi(\g))$ is the $SL(n,\R)$-orbit at $\g$.  
\een  
\eprop

\bcor 
Let $u\in C^\infty(S^1, V_n)$,  and $g:S^1\to GL(n,\R)$ such that $g^{-1}g_x= b+u$, where $b=\sum_{i=1}^{n-1} e_{i+1,i}$. Let $v\in C^\infty(S^1,V_n^t)$. If $(J_2)_u(v)=0$, then there is a constant $c_0\in sl(n,\R)$ such that $v= \pi_0(g^{-1}c_0g)$, where $\pi_0:sl(n,\R)\to V_n^t$ is the canonical projection.
\ecor

\begin{proof}
Let $\g$ denote the first column of $g$. The equation $g_x= g(b+u)$ implies that $g=(\g, \g_x, \cdots, \g_x^{(n-1)})$. By Corollary \ref{sy}, there exist $\d \g \in T\cm_n(\R)_\g$ such that $P_u(v)=g^{-1}\d g$, where $\d g=(\d \g, \cdots, (\d \g)_x^{(n-1)})$. So 
$$(J_2)_u(v)=[\p_x+b+u, g^{-1}\d g]=\rd \Psi_\g(\d \g)=0.$$ 
By Proposition \ref{fr}, there exists content $c_0 \in sl(n, \R)$ such that $\d \g=c_0 \g$. Therefore
$$\d g=(c_0\g, \cdots, (c_0 \g)_x^{(n-1)})=c_0(\g, \ldots, \g_x^{(n-1)})=c_0g.$$
This proves that $g^{-1}\d g=g^{-1} c_0 g$. Since $v=\pi_0(P_u(v))$, $v=\pi_0(g^{-1}c_0 g)$. 
\end{proof}

It is known that $\{\, , \}_1$ and $\{\, , \}_2$ are compatible (cf. \cite{DS84}), i.e. $c_1\{\, , \}_1+ c_2\{\, ,\}_2$ is a Poisson structure on $C^\infty(S^1, V_n)$ for all real constants $c_1, c_2$. It is standard in the literature (cf. \cite{Mag78}, \cite{TW1}) that we can use these two compatible Poisson structures to generate a sequence of Poisson structures: 
$$\{F_1, F_2\}_j(u)= \li (J_j)_u(\K F_1(u)), \K F_2(u)\ri,$$
where 
\beq\label{ix}
J_j=J_2(J_1^{-1}J_2)^{j-2}.
\eeq
Moreover, the $(nk+j)$-th $\an1$-KdV flow is
$$
u_{t_{nk+j}}= (J_{k+2})_u(\K H_j(u))= J_2(J_1^{-1}J_2)^k(\K H_j(u)),
$$
where $H_j$ is the functional on $C^{\infty}(S^1, V_n)$ defined by \eqref{sz}. 

\bs
\section{Bi-Hamiltonian Structure for central affine curve flows}\label{go2}

The pull back $\{\, ,\}_j^{\wedge}$ of the Poisson structure $\{\, , \}_j$ to $\calM_n(S^1)$ via the central affine curvature map $\Psi$ is for functions of the form $F\circ \Psi$, where $F$ is a functional on $C^\infty(S^1,V_n)$. In other words, $\{\, ,\}^\wedge_j$ is defined by
$$\{F_1\circ\Psi, F_2\circ \Psi\}^\wedge_j= \{F_1, F_2\}_j\circ \Psi$$
for functionals $F_1, F_2$ on $C^\infty(S^1, V_n)$.  
 Let $\hat w_j$ be the $2$-form defined by $\{\, ,\}_j^\wedge$, i.e., 
\beq\label{tb}
(\hat{w}_j)_{\g}(X_1 (\g), X_2 (\g))=\{\hat F_1, \hat F_2\}_j^\wedge(\g),
\eeq
where $X_i (\g)$ is the Hamiltonian vector field for $\hat F_i=F_i \circ \Psi$ with respect to $\{\, , \}^\wedge_j$, $i=1, 2$. In this section, we
\ben
\item show that $\hat w_2$ and $\hat w_3$ are the pull backs of symplectic forms on certain co-Adjoint orbits.  
\item prove that the 2-forms $\hat w_2$ and $\hat w_3$ induce symplectic forms on the orbit spaces $\frac{\calM_n(S^1)}{SL(n,\R)}$ and $\frac{\calM_n(S^1)}{SL(n,\R)\times \R^{n-1}}$ respectively, where $SL(n,\R)$ acts on $\calM_n(S^1)$ by $c\cdot \g= c\g$ and the $\R^{n-1}$-action on $\calM_n(S^1)$ is  generated by the first $(n-1)$ central affine curve flows \eqref{ncfj}.
\een

\bprop\label{iw}
\ni Fix $j$. Let $H$ be a functional on $C^\infty(S^1, V_n)$, $\Psi$ the central affine curvature map, and $\d$ the Hamiltonian vector field  for $\hat H= H\circ \Psi$ with respect to $\{\, ,\}_j^\wedge$.  Then for $\g\in \cm_n(S^1)$ we have
$$[\p_x+ b+u, g^{-1}\d g]= (J_j)_u(\K H(u)),$$
 where $g$ and $u$ are the central affine moving frame and central affine curvature along $\g$ respectively and $\d g=(\d\g, \ldots, (\d \g)_x^{(n-1)})$. 
\eprop

\begin{proof}  Since $\{\, ,\}_j^\wedge$ is the pull back of $\{\, , \}_j$, we have  
$$\rd \Psi_\g(\d \g)= (J_j)_u(\K H(u)).$$
By Proposition \ref{fp},  $\rd\Psi_\g(\d\g)= [\p_x+b+u, g^{-1}\d g]$.
\end{proof}

\bcor\label{hg} Let $H_j$ be functionals on $C^\infty(S^1,V_n)$ defined by \eqref{sz}, and $\hat H_j=H_j \circ \Psi$. Then the $j$-th central affine curve flow \eqref{ncfj} is the Hamiltonian equation for $\hat H_j$ and $\hat H_{n+j}$ with respect to $\{\, , \}_2^\wedge$ and $\{\, , \}_1^\wedge$ respectively. 
\ecor

\bcor Let $1\leq j\leq n-1$ and $k\geq 0$. Then
the $(nk+j)$-th central affine curve flow on $\calM_n(S^1)$ is 
$$\g_{t_{nk+j}} = gZ_{nk+j, 0}(u)e_1 = g(P_u((J_1^{-1}J_2)_u^k(\K H_j(u))))e_1,$$
where $g$ and $u$ are the central affine moving frame and curvature of $\g$ respectively. 
\ecor

\bprop\label{jn} Let $F_i$ be functionals on $C^\infty(S^1, V_n)$, and $\d_i$ the Hamiltonian vector field for $\hat F_i=F_i\circ \Psi$ with respect to $\{\,, \}^\wedge_j$ for $i=1,2$. Then
$$
\{\hat F_1, \hat F_2\}_j^\wedge(\g)= -\li g^{-1}\d_1 g, (J_2J_j^{-1}J_2)_u(\pi_0(g^{-1}\d_2 g))\ri,
$$
where $g$ is the central affine moving frame along $\g$, $u=\Psi(\g)$, and $\d_i g= (\d_i \g, (\d_i \g)_x, \ldots, (\d_i \g)_x^{(n-1)})$.
\eprop

\begin{proof} By Proposition \ref{fp}, $\rd \Psi_\g(\d \g)= (J_2)_u(\pi_0(g^{-1}\d g))$.  Since $\rd \Psi(\d_i\g)= (J_j)_u(\K F_i(u))$ for $i=1,2$, we get $\K F_2(u)= (J_j^{-1}J_2)_u(\pi_0(g^{-1}\d_2 g)$.  So we have $(J_j)_u(\K F_2(u))= J_2(\pi_0(g^{-1}\d_i g))$. Compute directly to get
\begin{align*}
\{\hat F_1,\hat F_2\}_j^\wedge (\g) & = \{F_1, F_2\}_j (u)=\li (J_j)_u(\K F_1(u)), \K F_2(u)\ri\\
&= \li J_2(\pi_0(g^{-1}\d_1 g), (J_j^{-1}J_2)_u(\pi_0(g^{-1}\d_2 g\ri\\
&=-\li \pi_0(g^{-1}\d_1 g), (J_2 J_j^{-1}J_2)_u(\pi_0(g^{-1}\d_2g))\ri\\
&=-\li g^{-1}\d_1 g,  (J_2 J_j^{-1}J_2)_u(\pi_0(g^{-1}\d_2g))\ri.
\end{align*}
\end{proof}

\bprop \label{is} Let $\hat w_j$ be the $2$-form on $\calM_n(S^1)$ defined by \eqref{tb},  i.e., 
$$
(\hat w_j)_\g(\d_1\g, \d_2 \g)= -\li g^{-1}\d_1 g, (J_2J_j^{-1}J_2)_u(\pi_0(g^{-1}\d_2 g))\ri,
$$
where $\d_i g=(\d_i \g, \cdots, (\d_i \g)_x^{(n-1)})$ for $i=1, 2$, and $g$ is the central affine moving frame along $\g$.
 Then
\begin{align}
(\hat w_2)_\g(\d_1\g, \d_2\g)&= \li [\p_x+b+u, g^{-1}\d_1g], g^{-1}\d_2 g\ri, \label{is2}\\
(\hat w_3)_\g(\d_1\g, \d_2\g)&=\li [e_{1n}, g^{-1}\d_1 g], g^{-1}\d_2 g\ri. \label{is3}
\end{align}
\eprop

Next we write down $\hat w_2(X,Y)$ and $\hat w_3(X,Y)$ in terms of determinants involving $X, Y\in T\calM_n(S^1)$ and derivatives of $X$ and $Y$: 
\bthm\label{hh} Let $X, Y$ be tangent vectors of $\calM_n(S^1)$ at $\g$.  Then
\begin{align*}
&(\hat w_2)_\g(X,Y)= -\sum_{i=1}^{n-1} \oint \det(\g, \ldots, \g_x^{i-2}, X_x^{(n)}, \g_x^{(i)}, \ldots, Y_x^{(i-1)})\rd x\\ 
&\quad +\sum_{i, j=1}^{n-1} \oint u_j\det(\g, \ldots, \g_x^{(i-1)}, X_x^{(j-1)}, \g_x^{(i)}, \ldots, Y^{(i-1)}_x)\rd x,\\
&(\hat w_3)_\g(X, Y)=- \sum_{i=1}^{n-1} \oint \det(\g, \ldots, \g_x^{(i-2)}, X, \g_x^{(i)}, \ldots, Y^{(i-1)}_x)\rd x\\
&\qquad -\sum_{i=1}^{n-1}\oint \det(\g, \ldots, \g_x^{(i-2)}, X_x^{(i-1)}, \g_x^{(i)}, \ldots, Y)\rd x. 
\end{align*}
\ethm 

\begin{proof}
Let $\g, g, u, \d_i\g, \d_i g$ be as in Proposition \ref{jn}, 
$$C=(C_{ij})= g^{-1}\d_1 g, \quad D=(D_{ij})= g^{-1}\d_2 g,$$
 and $C_i, D_i$ the $i$-th column of $C$ and $D$ respectively. By Theorem \ref{jw}, we can express $C_i$'s as differential polynomials in $C_1$. Similarly, $D_i$'s can be expressed as differential polynomials in $D_1$.  Moreover, $C_i=(C_{1i}, \ldots, C_{ni})$ is the coordinate of $(\d_1\g)_x^{(i-1)}$ with respect to the frame $g=(\g, \ldots, \g_x^{(n-1)})$, i.e., $(\d_1\g)_x^{(i-1)}= \sum_{k=1}^n C_{ki} \g_x^{(k-1)}$. 

Recall that if $Y= \sum_{i=1}^n y_i \g_x^{(i-1)}$, then 
$Y'= \sum_{i=1}^n (y_i'+ y_{i-1} + u_i y_n)\g_x^{(i-1)}$.  
Write $(\d_1\g)_x^{(n)} = \sum_{i=1}^n \xi_i (\d_1\g)^{(i-1)}_x$, 
Then $\xi=(\xi_1, \cdots, \xi_n)^t=C_n'+ (b+u) C_n$. 
By Proposition \ref{is}, we have
\begin{align*}
(\hat w_3)_\g(\d_1 \g, \d_2 \g) & =  \oint \sum_{i=1}^n C_{ni}D_{i1}- C_{i1} D_{ni} \rd x,\\
 (\hat w_2)_\g(\d_1\g, \d_2 \g) & = \oint \sum_{i=1}^n (C(b+u))_{in} D_{ni}-(C_x+(b+u)C)_{in} D_{ni}  \rd x\\
& = \oint \sum_{i=1}^n \sum_{j=1}^{n-1} C_{ij} u_j D_{ni}-\sum_{i=1}^n \xi_i D_{ni} \rd x .
\end{align*}

We compute $(\hat w_3)_\g$ as follows:  Let $X=\d_1 \g$, $Y=\d_2 \g$, then
\begin{align*}
&(\hat w_3)_\g(X,Y)= \oint \sum_{i=1}^n C_{ni}D_{i1}- C_{i1} D_{ni} \rd x\\
& =\oint C_{nn} D_{n1}-C_{n1} D_{nn} + \sum_{i=1}^{n-1} C_{ni} D_{i1}-C_{i1} D_{ni} \rd x\\
&= \oint C_{n1}\sum_{i=1}^{n-1} D_{ii} -(\sum_{i=1}^{n-1} C_{ii}) D_{n1} + \sum_{i=1}^{n-1}C_{ni} D_{i1}- C_{i1} D_{ni} \rd x.
\end{align*}
Note that $\det(\g, \ldots, \g_x^{(n-1)})=1$ and the $k$-th column of $C$ and $D$ are the coefficients of $X_x^{(k-1)}$ and $Y_x^{(k-1)}$ written as a linear combination of $\g, \ldots, \g_x^{(n-1)}$. So we have 
\beq\label{hs}
\det(\g, \g_x, \ldots, \g_x^{(i-2)}, X_x^{(k-1)}, \g_x^{(i)}, \ldots, Y_x^{(\ell-1)})= C_{ik} D_{n\ell} - C_{nk} D_{i\ell}.
\eeq
Substitute \eqref{hs} into the above formula for $w_\g(X,Y)$ to get the formula for $(\hat w_3)_\g$ as stated in the theorem.
 
Use  $\tr(C)=\tr(D)=0$ and \eqref{hs} to get the formula for $\hat w_2$. 
\end{proof}

\beg
For $n=2$, Theorem \ref{hh} gives
\begin{align*}
&(\hat w_2)_\g(X, Y)= -\oint \det(X', Y') + u_1\det(X,Y)\rd x, \\
&(\hat w_3)_\g(X,Y)=-2 \oint \det(X,Y)\rd x.
\end{align*}
These are the $2$ forms given in \cite{FK13} and \cite{UP95} respectively. 
\eeg

\beg
For  $n=3$, we get
$$(\hat w_3)_\g(X,Y)= -3\oint \det(X, \g', Y)\rd x,$$
which is the $2$ form given in \cite{CIM13}.   We also have 
\begin{align*}
(\hat w_2)_\g(X, Y)&=- \oint \det(X''', \g', Y) +\det(\g, X''', Y') \rd x \\
&\quad + \oint u_1(\det(X, \g', Y) +\det(\g, X, Y'))\rd x\\
&\quad  +\oint u_2(\det(X', \g', Y)+ \det(\g, X', Y'))\rd x.
\end{align*}
\eeg

\beg
For $n=4$, let $|v_1, v_2, v_3, v_4|=det(v_1, v_2, v_3, v_4)$.  Then we have
\begin{align*}
(\hat w_3)_\g(X, Y) &= -2\oint |X, \g', \g'', Y|+|\g, \g', X', Y'| \rd x,\\
(\hat w_2)_\g(X,Y) &=-\oint |X^{(4)}, \g', \g'', Y|+|\g, X^{(4)}, \g'', Y'|+|\g, \g', X^{(4)}, Y''|  \rd x\\
&\quad+ \oint u_3(|X'', \g', \g'', Y|+|\g, X'', \g'', Y'|+|\g, \g', X'', Y''|)\rd x \\
& \quad+ \oint u_2(|X', \g', \g'', Y|+|\g, X', \g'', Y'|+|\g, \g', X', Y''|) \rd x \\
& \quad+  \oint  u_1(|X, \g', \g'', Y|+|\g, X, \g'', Y'|+|\g, \g', X, Y''|) \rd x.
\end{align*}
\eeg

Next we prove that $\hat w_2$ and $\hat w_3$ are the pull backs of certain co-Adjoint orbit symplectic forms. Let $M$ be the co-Adjoint orbit of $G$ on the dual $\calG^*$ of the Lie algebra $\calG$ at $\ell_0\in \calG^*$. The orbit symplectic form on $M$ is defined by
$$\tau_\ell(\ti\xi(\ell), \ti\eta(\ell))= \ell([\xi,\eta]),$$ 
where $\ell\in M$, and $\ti \xi$ and $\ti\eta$ are infinitesimal vector fields corresponding to the co-Adjoint action generated by $\xi,\eta\in \calG$. 

We identify the Adjoint orbit of $C^\infty(S^1, SL(n,\R))$ on $C^\infty(\R, sl(n,\R))$ as the co-Adjoint orbit of $C^\infty(S^1,SL(n,\R))$ via the  non-degenerate  bi-linear form \eqref{qa}, i.e.,
$$\li \xi, \eta\ri =\oint \tr(\xi(x)\eta(x))\rd x.$$

\bthm\label{kg} Let $\calO_1$ denote the Adjoint orbit of $C^\infty(S^1,SL(n,\R))$ on $C^\infty(S^1, sl(n,\R))$ at the constant loop $e_{1n}$,  $\tau_1$ the orbit symplectic form on $\calO_1$, and $\fk_1$ the map from $\calM_n(S^1)$ to $\calO_1$ defined by $\fk_1(\g)= ge_{1n} g^{-1}$, where $g$ is the central affine moving frame along $\g$. Let $\hat w_3$ be as in \eqref{is3}.  Then  
$\fk_1^\ast\tau_1=\hat w_3$.
\ethm

\begin{proof}
Given $\xi\in C^\infty(S^1,sl(n,\R))$, a direct computation implies that the infinitesimal vector field is 
$$\ti \xi(ge_{1n}g^{-1})= [\xi, ge_{1n}g^{-1}].$$ 
So the orbit symplectic structure is
$$(\tau_1)_{ge_{1n}g^{-1}}([\xi, ge_{1n}g^{-1}], [\eta, ge_{1n}g^{-1}])= \li ge_{1n} g^{-1}, [\xi, \eta]\ri.$$
The differential of $\fk_1$ at $\g$ is 
$$\rd( \fk_1)_\g(\d\g)= [\d g g^{-1}, ge_{1n}g^{-1}].$$
Hence 
\begin{align*}
&(\fk_1^*\tau_1)_\g(\d_1 \g, \d_2 \g)= (\tau_1)_{ge_{1n}g^{-1}}([(\d_1 g)g^{-1}, ge_{1n}g^{-1}], [(\d_2 g)g^{-1}, ge_{1n}g^{-1}])\\
&= \li ge_{1n}g^{-1}, [(\d_1 g)g^{-1}, (\d_2 g)g^{-1}]\ri= \li e_{1n}, [g^{-1}\d_1 g, g^{-1}\d_2 g]\ri,
\end{align*}
which is equal to $(\hat w_3)_\g(\d_1\g, \d_2\g)$. 
\end{proof}

Let $\R \p_x+C^{\infty}(S^1, sl(n, \R))$ denote the Lie algebra with bracket defined by $$[r_1 \p_x+u, r_2 \p_x+v]=r_1v_x-r_2 u_x+[u, v], \quad r_1, r_2 \in \R.$$ 
It is known (cf. \cite{PS86}, \cite{Ter97}) that the dual of the central extension of the loop algebra $C^\infty(S^1,sl(n,\R))$ defined by the $2$-cocycle 
$$\rho(\xi,\eta)=\oint \tr(\xi_x(x)\eta(x))\rd x$$
can be identified as the Lie algebra $\R \p_x+ C^\infty(S^1,sl(n,\R))$. The co-Adjoint action corresponds to the gauge action,
$$g\cdot (\p_x+ u)=g(\p_x+u) g^{-1}= \p_x + gug^{-1} - g_xg^{-1}.$$

\bthm\label{kh} Let $\calO_2$ denote the gauge orbit of $C^\infty(S^1, SL(n,\R))$ at $\p_x$, $\tau_2$ the orbit symplectic form on $\calO_2$, and $\fk_2:\calM_n(S^1)\to \calO_2$  the map defined by 
$\fk_2(\g)= g^{-1}g_x$, where $g$ is the central affine moving frame along $\g$. Let $\hat w_2$ be as in \eqref{is2}.  Then $\fk_2^*(\tau_2)= \hat w_2$. 
\ethm

\begin{proof}
The infinitesimal vector field on $\calO_2$ given  by the gauge action for $\xi\in C^\infty(S^1,sl(n,\R))$ is
$\ti \xi(\p_x+ v)=- [\p_x+v, \xi]$.  
Note that 
$$\fk_2(\g)= \p_x+ g^{-1}g_x= \p_x+ b+ u= \p_x+ b+ \Psi(\g).$$ 
  By Proposition \ref{fp},  
$$\rd (\fk_2)_\g(\d \g)= \rd \Psi_\g(\d\g)= [\p_x+b+ u, g^{-1}\d g].$$  Then 
\begin{align*}
&(\fk_2^*\tau_2)_\g(\d_1 \g, \d_2\g)= (\tau_2)_{\p_x+ b+u}(\rd \fk_2(\d_1\g), \rd \fk_2(\d_2\g))\\
&=(\tau_2)_{\p_x+b+u}([\p_x+b+u, g^{-1}\d_1 g], [\p_x+ b+u, g^{-1}\d_2 g]) \\
&= \li [\p_x+b+u, g^{-1}\d_1 g], g^{-1}\d_2 g\ri,
\end{align*}
which is equal to $(\hat w_2)_\g(\d_1 \g, \d_2\g)$.
\end{proof}

Recall that a {\it weak symplectic form\/} on $M$ is a closed $2$-form on $M$ such that $w_x(v_1, v_2)=0$ for all $v_2\in TM_x$ implies that $v_1=0$. If $M$ is of finite dimension, then a weak symplectic form is symplectic. When $M$ is of infinite dimension, a weak symplectic form need not to be non-degenerated, but we can still have the Hamiltonian theory (cf. \cite{AM}). Below we show that $\hat w_2$ and $\hat w_3$ induce weak symplectic forms on the orbit spaces $\calM_n(S^1)/SL(n,\R)$ and $\calM_n(S^1)/(SL(n, \R)\times\R^{n-1})$ respectively.   

\bthm
The $2$-form $\hat w_2$  induces a weak symplectic form on the orbit space $\calM_n(S^1)/SL(n,\R)$. 
\ethm

\begin{proof} By Theorem \ref{kh}, $\hat w_2$ is a closed $2$-form. 
It follows from \eqref{is2} that $(\hat w_2)_\g(\d_1\g, \d_2\g)=0$ for all $\d_2\g$ if and only if $[\p_x+b+u, g^{-1}\d_1g]=0$, where $g$ is the central affine moving frame along $\g$ and $u=\Psi(\g)$.  By Proposition \ref{fp}, 
$$\rd \Psi_\g(\d_1\g)=[\p_x+b+u, g^{-1}\d_1 g]= (J_2)_u(\pi_0(g^{-1}\d_1 g)=0.$$ The theorem follows from Proposition \ref{fr}.
\end{proof}

We consider the $\R^{n-1}$-action on $\calM_n(S^1)$ generated by the first $(n-1)$ central affine curve flow \eqref{ncfj}. Since the central affine curve flows commute with the $SL(n,\R)$-action, the product group $SL(n,\R)\times\R^{n-1}$ acts on $\calM_n(S^1)$.  

\bthm The $2$-form
 $\hat w_3$ induces a weak symplectic form on the space $\calM_n(S^1)/(SL(n,\R)\times \R^{n-1})$.
 \ethm

\begin{proof} By Theorem \ref{kg}, $\hat w_3$ is a closed $2$-form.  The formula of $\hat w_3$ implies that $(\hat w_3)_\g(\d_1 \g, \d_2 \g)=0$ for all $\d_2 \g$ in $T\calM(S^1)_\g$ if and only if 
$$(J_1)_u(\pi_0(g^{-1}\d_1 g))=0.$$ So $\pi_0(g^{-1}\d_1g)$ lies in the kernel of $(J_1)_u$. It follows from Theorem \ref{jr} that $\xi$ lies in the span of $\K H_1(u), \ldots, \K H_{n-1}(u)$, where $H_i$'s are the Hamiltonians defined by \eqref{sz}.  

Let $X_1, \ldots, X_{n-1}$ denote  the vector fields that generated the first $(n-1)$ central affine curve flows. Then $\rd \Psi_\g(X_i)=(J_2)_u(\K H_i(u))$.  
If $\pi_0(g^{-1}\d_1 g)= \K H_i(u)$ for some $1\leq i\leq n-1$, then $\rd\Psi(\d_1\g)= (J_2)_u(\K H_i(u))$. So $\d_1\g\in X_i(\g) + \Ker(\rd \Psi_\g)$.  By Proposition \ref{fr}, $\Ker(\rd\Psi_\g)=$ the tangent space of the $SL(n,\R)$-orbit at $\g$. This finishes the proof. 
\end{proof}

\bs

\end{document}